\documentclass[hidelinks,onefignum,onetabnum]{siamart251104}

\usepackage[utf8]{inputenc} 
\usepackage[T1]{fontenc}    
\usepackage{xcolor}
\usepackage{hyperref}       
\usepackage{url}            
\usepackage{booktabs}       
\usepackage{nicefrac}       
\usepackage{microtype}      
\usepackage{graphicx}
\usepackage{amsmath,amssymb}
\usepackage{enumitem}
\usepackage{bbm}
\usepackage{dsfont}
\usepackage{mathrsfs}
\usepackage{float}
\usepackage[sort]{cite}

\usepackage{algorithm}
\usepackage{algpseudocode}

\usepackage{subcaption}
\usepackage[labelformat=simple]{subcaption}

\numberwithin{equation}{section}

\newtheoremstyle{nonitalic}
{3pt}
{3pt}
{\normalfont}
{}
{\bfseries}
{}
{}

\theoremstyle{nonitalic}

\newcounter{hypothesis}

\DeclareMathOperator{\sinc}{sinc}
\DeclareMathAlphabet\mathbfcal{OMS}{cmsy}{b}{n}

\DeclareMathOperator{\supp}{supp}

\newcommand{\as}[1]{\left\vert#1\right\vert}

\newcommand{\norm}[1]{\left\Vert#1\right\Vert}

\newcommand{\grad}{\nabla}

\newcommand{\dx}{{\rm d}x}
\newcommand{\dy}{{\rm d}y}

\headers{Stationary States for Nonlocal Fokker-Planck Equations}{J. A. Carrillo \& Y. Salmaniw \& A. L. Villares}

\author{Jos\'e A.~Carrillo\thanks{\,Mathematical Institute, University of Oxford, Woodstock Road, Oxford OX2 6GG, United Kingdom.
(\email{jose.carrillo@maths.ox.ac.uk})}
\and
Yurij Salmaniw\thanks{\,Department of Mathematics, Physics and Geology, Cape Breton University, Sydney, Nova Scotia, Canada.
(\email{yurij\_salmaniw@cbu.ca})}
\and
Antonio Le\'on Villares\thanks{\,Department of Engineering Science, University of Oxford, Parks Road, Oxford OX1 3PJ, United Kingdom.
(\email{antonio.leonvillares@stx.ox.ac.uk})
}
}

\title{Numerical stationary states for nonlocal Fokker-Planck equations via fixed points of consistency maps}

\usepackage{amsopn}

\ifpdf
\hypersetup{
pdftitle={Numerical stationary states for nonlocal Fokker-Planck equations via fixed points of consistency maps},
pdfsubject={},
pdfauthor={J. A. Carrillo \& Y. Salmaniw \& A. L. Villares},
pdfkeywords={},
}
\fi

\begin{document}

\maketitle

\begin{abstract}
We propose a fixed-point-based numerical framework for computing stationary states of nonlocal Fokker-Planck-type equations. Instead of discretising the differential operators directly, we reformulate the stationary problem as a nonlinear fixed-point map built from the original PDE and its nonlocal interaction terms, and solve the resulting finite-dimensional problem with a matrix-free Newton-Krylov method. We compare implementations using the analytic Frechet derivative of this map with a simple central-difference approximation. Because the method does not rely on time evolution, it is agnostic to dynamical stability and can detect both stable and unstable stationary states. Its accuracy is determined mainly by the numerical treatment of convolutions and quadrature, rather than by differentiation stencils. We apply the approach to three model problems with linear diffusion, use existing analytical results to verify the outputs, and reproduce known bifurcation diagrams, as well as new bifurcation behaviour not previously observed in this kind of problem.
\end{abstract}

\section{Introduction}
In this paper, we seek to numerically identify the stationary states of problems of the form
\begin{align}\label{eq:time_pde_general}
\begin{cases}
        \partial_t u = \grad \cdot ( u \grad \left ( H^\prime [u] + \Phi[u](x) + V(x) \right) ), \quad x \in \Omega, \quad t > 0, \\
    u(x,0) = u_0(x) \geq 0 ,
\end{cases}
\end{align}
where $u(x,t) \geq 0$ is the unknown function (e.g., a density or probability density), $u_0(x)$ is the initial (probability) density, $H(u)$ is a density of internal energy, and $V(x)$ is a given potential function. The function $\Phi[u](x)$ encodes nonlocal, typically nonlinear, interactions, whose form varies depending on the application. $\Omega$ is a spatial domain that may be the whole space, a compact manifold without boundary, or some combination; the typical difficulty that arises is in defining the nonlocal quantity $\Phi[u](x)$ near a physical boundary $\partial \Omega$. Therefore, it is most common to study such equations on a domain without a boundary.

Equations of the form \eqref{eq:time_pde_general} appear in many contexts, including fluid dynamics \cite{paterson1981first}, phase separation in materials science \cite{HP06,MR1453735}, cell-cell adhesion in biology \cite{MR3783102,MR3948738,ButtenschoenHillen2021,Buttenschoen2018SpaceJump,falco2022local}, crowd dynamics \cite{muntean2014collective,MR1698215}, and biological aggregations \cite{MR2257718}. Their ubiquity can be explained by their direct connection with discrete dynamics described by stochastic differential equations \cite{Carrillo2020LongTime,pareschi2013interacting,MR3858403}, random walk processes \cite{potts2016territorial, pottslewis2016}, and singular limits connecting them with higher-order PDEs \cite{elbar2023limit, elbar-skrzeczkowski, MR4911888}. Here, the nonlocal component can model phenomena such as the ability to detect some stimulus far away from the point of detection. In cell biology, for example, this manifests as protrusions called \textit{filopodia} which allow cells to probe their environment, influencing their movement at the population level \cite{carrillo2019aggregation}. Other phenomenological models incorporate the ability of larger animals to detect scent marks of other competing populations or to remember spatial locations of previous species-species interactions \cite{pottslewis2019}, promoting the formation of territories \cite{potts2016territorial}; in other cases, systems of interacting populations have been used to demonstrate the emergence of so-called \textit{run-and-chase} dynamics \cite{MR3783102,Painter2024a, jewell2025chaseandrunchiralitynonlocalmodels}.

Motivated by their ubiquity and potential for broad application to the sciences, a comprehensive description of their long-time dynamics is of particular interest. From an analytical point of view, we can understand some of this behaviour using linearisation techniques, bifurcation theory, asymptotic analysis, and related tools from functional analysis and PDEs. In general, it is quite challenging to analytically determine quantitative solution behaviour in the long-time limit. Therefore, we combine these functional analytical techniques with a numerical scheme to directly solve for the stationary states of model \eqref{eq:time_pde_general} in a general and robust manner.

Finite-volume schemes preserving both the positivity and the dissipation of the free energy at the fully discrete level in time have been proposed in the literature, e.g., \cite{MR3023729,MR3372289,B.C.H2020,MR4605931}. They are among the salient numerical approaches to deal with evolutionary equations of the form \eqref{eq:time_pde_general}. These numerical methods combine techniques from systems of hyperbolic conservation laws, such as upwinding, with a careful discretisation of the velocity fields in order to capture the desired dissipative structure. Even though they are excellent numerical methods for studying the evolution of these problems, they have inherent limitations for computing stationary solutions. Any time-dependent solver will only detect those states that are dynamically stable, locally or otherwise. Therefore, bifurcation branches are not feasible to find with these solvers, and, moreover, stabilisation in time might be tricky if metastability is present in the system \cite{carrillo2024well}, due to different time scales in their evolution. It is also unclear how advantageous it is to solve for stationary states directly by iterative Newton-type methods compared to the time it takes for time-dependent solvers to reach equilibrium.

Our numerical approach does not require discretising spatial or temporal derivatives due to the structural properties of equations of the form \eqref{eq:time_pde_general}. Therefore, stability and accuracy properties of a discretisation scheme (e.g., CFL conditions for explicit solvers or Newton iterations for implicit solvers) are not a concern in the classical sense. 

\subsection{Existing tools \& methods}\label{sec:existing_tools_methods}

A standard workhorse for solving nonlinear algebraic systems of the form $F(U) = 0$ is a Newton-Krylov approach, comprising at least two key steps: the ``outer'' Newton portion, which solves the original nonlinear problem, and the ``inner'' Krylov portion, which solves the linear system appearing at every step of the Newton scheme. More precisely, denoting by $J$ the Fr{\'e}chet derivative of $F$, the outer stage updates
\(U_{k+1}=U_k+\delta U_k\) by solving the inner problem
$$
J(U_k)\,\delta U_k = -F(U_k),
$$
typically via Krylov-subspace approximation methods (e.g.\ GMRES \cite{gmres}), rather than a direct factorisation. This approach has been effectively implemented to study, for example, pattern formation in classical reaction-diffusion systems \cite{MR2929933}, general dissipative dynamical systems \cite{UmbriaNet2016}, and other kinds of nonlocal PDE through, e.g., fractional diffusion \cite{MR4189291} or in applications to neuroscience \cite{MR3164128}. This ``outer/inner'' viewpoint is emphasised in the Jacobian-free Newton-Krylov (JFNK) literature: the outer Newton mechanism is largely generic, while efficiency is dominated by the inner linear solver through preconditioning and by how Jacobian information is accessed \cite{MR2030471,MR3826508}.

Assuming that reduction to finite dimension has been handled effectively for the outer stage, the discretised stationary version of problem \eqref{eq:time_pde_general} falls squarely within this paradigm. In the PDE setting, the nonlinear map $F$ obtained inherits the structure of the original PDE; consequently, its Jacobian $J$ inherits properties of the linearised differential operator together with any lower-order, nonlinear, or nonlocal couplings. Therefore, the effective implementation of the \textit{inner} stage of the Newton-Krylov pipeline is highly problem-dependent, even within the narrow class of \eqref{eq:time_pde_general}. Even for problems with constant linear diffusion, one must carefully design preconditioners to allow for efficient scaling \cite{MR2929933}; as one moves away from such cases, we require substantially different solver and preconditioner choices \cite{UmbriaNet2016, MR3283356, MR4189291}. In general, there is no universal preconditioner: good performance typically depends on exploiting the specific operator structure, a central theme of modern preconditioning surveys \cite{MR4175150}.

Concerning the outer stage, whilst Newton's method is relatively simple, applying the pipeline to \eqref{eq:time_pde_general} directly requires discretising the PDE appropriately. Deriving a general numerical scheme that effectively discretises \eqref{eq:time_pde_general}, along with efficient preconditioning choices, does not appear in the literature and remains a challenge.

These challenges with the inner and outer stages are intrinsic to PDEs in general, and therefore existing software packages for numerical continuation of bifurcation branches are not a robust approach: given the generality of problems falling in the class of \eqref{eq:time_pde_general}, it is unlikely that they will handle generic nonlinear diffusion (through $H$) and general nonlocal terms (through $\Phi$).

To overcome some of these issues described above, we use the analytical structure of \eqref{eq:time_pde_general} to recast the stationary problem as an equivalent fixed-point equation, reducing some of the challenges related to the outer stage: discretisation will only require quadrature and convolutions, so that high-order accuracy is governed primarily by the underlying integration rules and FFT-based convolution, rather than differentiation schemes. This also reduces challenges with the inner stage: the Fr{\'e}chet derivative of the fixed-point map can often be written in closed form, providing guaranteed accuracy of the Jacobian at the inner stage, and reducing the computational cost when compared to a finite-difference approximation. 

Ultimately, recasting the stationary problem as a fixed-point map fundamentally changes the linear system arising at the inner stage of the Newton-Krylov approach. Our linearised operator takes into account that we have effectively inverted the differential operator a priori: our Fr{\'e}chet derivative is, roughly speaking, the identity plus a compact-like operator (see the Fr{\'e}chet derivative  \eqref{eq:mckean_vlasov_frechet_der} for the McKean-Vlasov case), and hence the resulting Jacobian is well-conditioned\footnote{It is not universally so, however; indeed, the Jacobian can become ill-conditioned, for example, near a bifurcation point.}. In particular, our approach appears to be highly effective at recovering stationary states for a variety of problems without requiring any preconditioning. We now describe our approach in more detail.

\subsection{Description of the approach}\label{sec:description_of_approach}

For now, we aim to provide a heuristic understanding of the setup, saving additional technical details for specific cases considered in later sections. Setting $\partial_t u = 0$ in \eqref{eq:time_pde_general}, one uses the conservation of mass and positivity of the solution to argue that any stationary state $u^*$ satisfying \eqref{eq:time_pde_general} is necessarily a fixed point of a nonlinear map $\mathcal{T} : X \mapsto X$ in the sense that
\begin{align}\label{eq:stationary_iff_fixedpoint}
    X \ni u^* = u^*(x) \text{ is a stationary state of \eqref{eq:time_pde_general}} \ \iff \ \mathcal{T}u^* - u^* = 0,
\end{align}
where $X$ is a Banach space dependent on the particular problem. For many problems taking the form of \eqref{eq:time_pde_general}, there is an associated free-energy functional that decays along solution trajectories\cite{carrillo2019aggregation}. This generally rules out persistent temporal oscillations in the long-time dynamics \cite{MR2763076}, and this fixed-point formulation will, in principle, capture \textit{all} possible limit points, stable or otherwise. Even in cases where there is no energy decay, one can still identify solutions through the fixed point; however, it is then possible that there are other kinds of solutions (e.g., periodic orbits) that are not captured by this fixed-point formulation.

To identify $\mathcal{T}$, one argues that with appropriately prescribed boundary conditions (e.g., no-flux or periodic) and mass constraints, any stationary state $u$ must satisfy
\begin{align}\label{eq:fixed_point_preliminary}
 H^\prime [u] + \Phi[u](x) + V(x) = c_0 \quad \text{ on }\quad  \supp (u) := \overline{\{ x \in \Omega : u(x) > 0 \}},
\end{align}
and $c_0 \in \mathbb{R}$ is a constant depending on $u$ which may change on different connected components of the support. For simplicity, we will always assume that $\supp (u) = \Omega$ so that there is only one constant $c_0$ to identify.

Assuming that $H^\prime : \mathbb{R}^+ \mapsto \mathbb{R}$ is invertible on its domain, we have that
\begin{align}
    u = ( H^\prime )^{-1} \left[ c_0(u) - \Phi[u](x) - V(x) \right],
\end{align}
where $( H^\prime )^{-1}[ \cdot ]$ denotes the functional inverse. Three revealing forms for $H$ include:
\begin{itemize}
    \item \textbf{Linear diffusion: \cite{MR1617171}} $H[u] = \sigma u ( \log u - 1 )$, $\sigma > 0$. Then, $H^\prime [u] = \sigma \log u$, and $(H^\prime)^{-1}[u] = \exp( {u/\sigma} )$.
    \item \textbf{Porous medium/fast diffusion: \cite{MR1260981}} $H[u] = \nu u^m$ for constants $\nu > 0$ and $m > 1$ (porous media) or $m \in (0,1)$ (fast diffusion). Then, $H^\prime [u] = \nu m u^{m-1}$, and $(H^\prime)^{-1} [u] = \left( u / \nu m \right)^{1/(m-1)}$.
    \item \textbf{Fermi-Dirac/exclusion/volume filling: \cite{MR2426970}} $H[u] = u\log u + (1-u)\log(1-u)$. Then, $H^\prime[u]=\log(\tfrac{u}{1-u})$ and $(H^\prime)^{-1} [u]=1/(1+\exp({-u}))$.
    \item \textbf{Bose-Einstein: \cite{MR3485127}}
$H[u]=u\log u-(1+u)\log(1+u)$. Then, $H'[u]=\log ( u / (1+u) )$ and $(H^\prime)^{-1} [u]=-1/(1-\exp({-u}))$. 
\end{itemize}
In general, we can identify the map $\mathcal{T}$ to be
\begin{align}
    \mathcal{T}u := ( H^\prime )^{-1} \left[ c_0(u) - \Phi[u](x) - V(x) \right].
\end{align}
However, we note that it is nontrivial to establish the equivalence property \eqref{eq:stationary_iff_fixedpoint} for general nonlinearities $H[u]$. More precisely, one always has the forward implication through construction (a stationary solution is a fixed point); the reverse direction (a fixed point is a stationary solution) requires further argumentation. For some cases (e.g., linear diffusion),  \eqref{eq:stationary_iff_fixedpoint} has been established \cite[Proposition 2.4]{Carrillo2020LongTime}. Hereafter, we take \eqref{eq:stationary_iff_fixedpoint} for granted, leaving open an interesting area of future analytical study.

Once $\mathcal{T}$ has been identified, we can translate the challenge of numerically solving for stationary states through the discretisation of the PDE into the relatively simpler problem of numerically finding roots of $\mathcal{T} - I$ by applying well-known optimisation techniques to its discretised form. In practice, a standard Newton method is not feasible due to the computational cost of constructing and inverting the Jacobian, so we opt for a Jacobian-free Newton-Krylov method \cite{MR2030471}. We can further improve the efficiency and stability of this approach by exploiting the analytical Fr{\'e}chet derivative of $\mathcal{T}$, at least for cases when the analytical computation is feasible; in this work, we choose instances such that this is possible so that we may compare analytical versus discrete implementations of the Fr{\'e}chet derivative.  

\subsection{Main contributions}

Our goals can be broadly summarised as follows. First, we develop and test a fixed-point-based numerical framework for computing stationary states of three representative nonlocal Fokker-Planck models with linear diffusion, each of which falls within the class of equations described by \eqref{eq:time_pde_general}. These examples, namely the McKean-Vlasov equation, a Cucker-Smale flocking model, and a neural Fokker-Planck equation, are chosen to demonstrate the robustness with respect to spatial domains, boundary conditions, and the nonlocality structure. Second, we validate the numerical output against existing analytical information: for the McKean-Vlasov case, this includes explicit stationary solutions and local bifurcation theory for special cases; for the other models, we use a combination of numerical reference solutions and local bifurcation theory results. Third, we compare the performance of implementation based on using the analytical Fr{\'e}chet derivative of the nonlinear map versus a central-difference approximation of the Fr{\'e}chet derivative. Finally, we use this framework to construct bifurcation diagrams for these models, both as a consistency check against known results and as a tool to reveal an extremely rich world of bifurcation behaviours and solution structures, at least for these kinds of equations.

These goals are realised in the remainder of the manuscript as follows. In Section \ref{sec:case-studies} we introduce the three model problems and, for each, derive the corresponding nonlinear map $\mathcal{T}$ and its relation to the original PDE, thereby illustrating how the fixed-point framework accommodates different domains and boundary conditions. Section \ref{sec:mckean-vlasov-eqn} then focuses on the McKean–Vlasov example in detail: we state sufficient conditions under which stationary states coincide with fixed points, compute the analytic Fr{\'e}chet derivative of $\mathcal{T}$, and describe the numerical implementation of the fixed-point solver, including the discretisation of the nonlocal terms. The main numerical results are then collected in Section \ref{sec:results}. In Section \ref{sec:kuramoto_model_1}, we validate the method for the McKean–Vlasov equation against analytical stationary solutions and known bifurcation diagrams for the Kuramoto model of synchronised oscillators, and compare the use of analytic versus finite-difference Fr{\'e}chet derivatives. Section \ref{sec:two_mode_mckean_vlasov} presents further bifurcation results, including super- and subcritical branches predicted by local theory, and new behaviour arising at a point where the linearised operator has a two-dimensional kernel. In Section \ref{sec:general_MVE} we explore additional McKean–Vlasov cases with interaction kernels commonly used in the literature and obtain bifurcation diagrams and stationary patterns that, to our knowledge, have not previously been reported. Section \ref{sec:init_iterate} discusses the influence of the initial iterate in the Newton scheme, while Section \ref{sec:twodim_MVE} illustrates that the approach extends naturally to two spatial dimensions by displaying stationary profiles for the McKean–Vlasov equation. Finally, in Section \ref{sec:CS_NFP_results} we present a brief validation for the Cucker–Smale and neural Fokker–Planck models, showing that the same fixed-point solver reproduces their expected stationary states and bifurcation structure.

We believe that this approach offers an excellent tool for exploring more deeply the behaviour of stationary solutions of equations of the form \eqref{eq:time_pde_general}. Our approach can be extended further beyond the scope of the present work, such as for cases with nonlinear diffusion (e.g., porous medium or volume-filling described above), or nonlocal systems of several interacting populations \cite{carrillo2025longtimebehaviourbifurcationanalysis}, with the caveat that additional structural requirements may need to be imposed to ensure the system cases fall within this fixed point framework.

\section{Nonlocal Fokker-Planck models: Case Studies}\label{sec:case-studies}

We choose three instances of model \eqref{eq:time_pde_general} to demonstrate the utility of this approach, namely, the \textit{McKean-Vlasov equation} \cite{Carrillo2020LongTime}, a \textit{Cucker-Smale flocking model} \cite{MR3541988}, and a \textit{neural Fokker-Planck equation} \cite{MR4491042, MR4575120}. These models are chosen for several reasons. First, they are compatible with the fixed-point framework outlined above. More precisely, we will choose linear diffusion at rate $\sigma > 0$ for all three cases:
$$
H [u] := \sigma u ( \log u - 1 ),
$$
so that the inverse of $H^\prime [u]$ is an exponential as described in Section \ref{sec:description_of_approach}. Second, these models have well-understood analytical properties that serve as robust verification and validation of the numerical output. Indeed, in some cases, we have access to an (implicit) analytical solution to test against; in all cases, we have local bifurcation theory results to verify the expected qualitative and quantitative solution behaviour. Finally, these models feature several different domains/boundary conditions, demonstrating that the approach is robust across different model configurations. We describe these three case studies in more detail now.

\subsection{The McKean-Vlasov Equation}\label{sec:dev-mckean-vlasov}

Fix the dimension $d=1$ and choose the domain $\Omega = \mathbb{T} := \mathbb{R}/L\mathbb{Z}$ to be the torus of side-length $L>0$, and fix
$$
\Phi[u](x) := \kappa \, W*u, \quad V(x) \equiv 0,
$$
where $W*\cdot$ denotes a (periodic) spatial convolution with a given sufficiently smooth and even interaction kernel $W$, and $\kappa \geq 0$ is an interaction strength parameter. This yields the stationary McKean-Vlasov equation on the torus \cite{Carrillo2020LongTime}
\begin{align}\label{eq:McKean_Vlasov_stationary}
        0 =  \partial_x \left( u \partial_x \left( \sigma \log u +  \kappa W * u \right) \right), \quad x \in \mathbb{T},
\end{align}
subject to some mass constraints. Solutions are then fixed points of the nonlinear map $\mathcal{T}:L^2(\mathbb{T}) \mapsto L^2(\mathbb{R})$ given by
\begin{align}\label{eq:operator_T_McKean}
    \mathcal{T}u :=  \frac{ \exp \left( - \sigma^{-1} \kappa W*u \right) }{Z(u)} , \quad Z(u) := \int_\mathbb{T} \exp \left( - \sigma^{-1} \kappa W*u \right) \dx . 
\end{align}
When the kernel $W$ is of a particular form, an implicit analytical solution is available to test against. We will also apply our approach to the dimension $d=2$ case of \eqref{eq:McKean_Vlasov_stationary}.

\subsection{A Cucker-Smale Flocking Model}\label{sec:dev-cucker-smale}
Fix $d = 1$, choose $\Omega = \mathbb{R}$, and fix
\begin{align}
    \Phi[u](x) := W*u, \quad V(x) := \alpha \left( \frac{\as{x}^4}{4} - \frac{\as{x}^2}{2}\right),
\end{align}
where we choose an interaction potential of confinement type, namely, $W(x) = \as{x}^2 / 2$. This yields the stationary problem 
\begin{align}\label{eq:cucker_smale_stationary}
    0 = \partial_x \left( u \, \partial_x \left( \sigma \log u + W*u + \alpha \left[ \tfrac{\as{x}^4}{4} - \tfrac{\as{x}^2}{2} \right] \right) \right), \quad x \in \mathbb{R},
\end{align}
subject to some mass constraints, with solutions of the form $\mathcal{T}: \mathbb{R} \mapsto L^1(\mathbb{R})$ given by
\begin{align}
    \mathcal{T} [\overline{u}](x) := &\frac{\exp \left( - \sigma^{-1} \left[ \alpha \frac{\as{x}^4}{4} + (1-\alpha) \frac{\as{x}^2}{2} - \overline{u} x \right] \right)}{Z(\overline{u}, \sigma)}, \nonumber \\
     Z(\overline{u}, \sigma) := &\int_{\mathbb{R}^d} \exp \left( - \sigma^{-1} \left[ \alpha \frac{\as{y}^4}{4} + (1-\alpha) \frac{\as{y}^2}{2} - \overline{u} y \right] \right) \dy,
\end{align}
where $\overline{u}$ is a mean velocity to be determined through the relation
$$
\overline{u} = \int_\mathbb{R} x \mathcal{T}[\overline{u}](x) \dx.
$$
As observed in \cite{MR3541988}, stationary states can be fully parametrised by identifying the velocity $\overline{u}$ according to the fixed point relation above. This yields an implicit analytical solution to test against.

\subsection{A Neural Fokker-Planck Model}\label{sec:dev-neural-fp}
Choose $\Omega = \mathbb{T} \times [0,\infty)$, $B>0$, and fix
\begin{align}
      \Phi[u](x) := x^2/2 - \int_0^x F\left( W * \overline{u} (y) + B  \right)\dy , \quad V \equiv 0,
\end{align}
where $W*\cdot$ is a spatial convolution, and $\overline{u}(x)$ denotes the first moment with respect to the second variable:
$$
\overline{u}(x) := \int_0^\infty y u(x,y) \dy .
$$
Here, $F(\cdot)$ describes the activation level, typically a saturating function like ReLU or the sigmoidal form \cite{MR4491042}, and $B$ is a constant external input signal. Then, $x \in \mathbb{T}$ represents the angular domain, and $u(x,y)$ describes the probability of observing the activity level $y \in [0,\infty)$ at a stationary state.

This yields the stationary problem for $(x,y) \in \mathbb{T} \times (0,\infty)$:
\begin{align}\label{eq:neural_FP_stationary}
    0 = \partial_x \left( u \partial_x \left( \sigma \log u + x^2/2 - \int_0^x F\left( W * \overline{u} (y) + B  \right)\dy\right) \right).
\end{align}

Different from the first two examples, this model is not of gradient-flow type. However, it is assumed that solutions are non-negative when equipped with no-flux boundary conditions at $y=0$:
$$
\left[ F( W* \overline{u}(x) + B ) u(x,y) - \sigma \partial_y u(x,y) \right]\bigr\vert_{y = 0} = 0,
$$
and sufficient decay as $y \to +\infty$ so that the mass is preserved. Denote $F_0(x) :=  F( W* \overline{u}(x) + B)$. Then, stationary solutions are fixed points of the map
\begin{align}
    \mathcal{T} u = \frac{\exp \left( - \tfrac{( y - F_0(x))^2}{2 \sigma} \right)}{Z (u,x)}, \quad Z(u,x) := L \int_0^\infty \exp \left( - \tfrac{( y - F_0(x))^2}{2 \sigma} \right) \dy
\end{align}
This normalises solutions so that $\int_0^\infty u(x,y) \dy = L^{-1}$ for each angle $x \in \mathbb{T}$, and so that upon integration over the angular domain, $u(x,y)$ is a probability density function. It is possible to characterise those stationary states which are homogeneous in the angular domain; under some additional technical criteria (see, e.g., \cite[Proposition 3.1]{MR4491042}), the solution obtained from the mappings above is the unique spatially homogeneous stationary solution. While not exactly of the form displayed in \eqref{eq:time_pde_general}, it is compatible with a fixed-point formulation and occurs on a semi-infinite domain with a no-flux boundary condition. A local bifurcation analysis carried out in \cite{MR4575120} serves as further validation of our numerical approach.

Together, these three case studies serve as excellent examples to demonstrate the accuracy, efficiency, and robustness of the approach as we shall exhibit in Section \ref{sec:results}. Before presenting our main results, we first develop in significant detail the numerical approach for the first case study: the McKean-Vlasov equation.

\section{The McKean-Vlasov Equation}\label{sec:mckean-vlasov-eqn}

\subsection{Recasting the problem as a pointwise evaluation map}

An essential part of our approach is the ability to recast the stationary PDE as a fixed-point problem. This way, rather than discretising the differential operator and projecting onto a finite-dimensional subspace, we work directly with this nonlinear map. This avoids some of the issues raised in Section \ref{sec:existing_tools_methods}, offering an operator-free approach in the traditional sense, and potentially providing scalability and accuracy improvements over more general-case bifurcation software packages. We explain this approach in detail for the McKean-Vlasov model \eqref{eq:McKean_Vlasov_stationary} as follows.

Equations of the form \eqref{eq:time_pde_general} are known to have a free energy whose time-derivative decreases along solution trajectories \cite{carrillo2019aggregation}. The McKean-Vlasov equation falls squarely into this paradigm, and solutions to the stationary problem \eqref{eq:McKean_Vlasov_stationary} necessarily satisfy the generic relation \eqref{eq:fixed_point_preliminary}. Moreover, primarily due to the presence of linear diffusion, any weak solution is guaranteed to be smooth and strictly positive in $\mathbb{T}$. Therefore, there is a single constant $c_0$ valid over the entire domain (i.e., the support of $u$ is the entire domain), and the constant is given by $Z(u)$ defined in \eqref{eq:operator_T_McKean}. Consequently, the fixed point formulation $\mathcal{T}u = u$ also defined in \eqref{eq:operator_T_McKean} is valid, and we have the equivalence introduced in \eqref{eq:stationary_iff_fixedpoint}. This is the nonlinear map we will use to numerically obtain the stationary solution(s) through identification of zeros of the functional $\mathcal{F}u := \mathcal{T}u - u$.

\subsection{Fr{\'e}chet Derivatives}

In order to apply Newton-type algorithms, one requires access to the Fr{\'e}chet derivative of the map $\mathcal{F}$ along with its inverse; this amounts to identifying the Fr{\'e}chet derivative of the nontrivial component $\mathcal{T}$, which we denote by $D \mathcal{T}[u](\phi)$. Here, $D\mathcal{T}[u](\phi)$ is understood as the Fr{\'e}chet derivative of $\mathcal{T}$, evaluated at $u$, in the direction $\phi$, i.e., the map $D\mathcal{T}[u](\phi)$ is a bounded linear functional from $L^2$ into itself, and for any mean-zero variation $\phi$ there holds
$$
\norm{\mathcal{T}( u + \phi) - \mathcal{T}u - D \mathcal{T}[u] (\phi)} = o ( \norm{\phi} ) ,
$$
where $\norm{\cdot}$ denotes the $L^2(\mathbb{T})$-norm.

In practice, we can approximate the inverse numerically through the GMRES algorithm \cite{gmres}, which only requires the ability to evaluate the discretisation of Fr{\'e}chet derivative $D\mathcal{T}[u](\phi)$ (i.e., the Jacobian of the discretised map) at some vector $\phi$. If we have access to the analytical Fr{\'e}chet derivative, this is straightforward; otherwise, a finite-difference scheme can be used to approximate the Fr{\'e}chet derivative numerically.

For the McKean-Vlasov equation, we have access to the analytical Fr{\'e}chet derivative, though we acknowledge that in general this computation can vary in difficulty, depending on the particular problem at hand. Therefore, we compare its direct use versus a finite-difference approximation without preconditioners. For the map $\mathcal{T}$ identified in \eqref{eq:operator_T_McKean}, its Fr{\'e}chet derivative $D \mathcal{T}$ at $u$ in the direction $\phi$ is given by
\begin{align}\label{eq:mckean_vlasov_frechet_der}
    D \mathcal{T}[u] (\phi) = - \sigma^{-1} \kappa \mathcal{T}u \left( W * \phi - \int_\mathbb{T} \mathcal{T}u \, W * \phi \, \dx \right),  
\end{align}
Therefore, the Fr{\'e}chet derivative of the map $F[u] := \mathcal{T}[u] - u$ is simply 
$DF[u](\phi) = D \mathcal{T}[u] (\phi) - I(\phi)$,
where $I$ is the identity map. For explicit computation of this quantity, along with higher-order Fr{\'e}chet derivatives, we refer readers to \cite[Appendix A.2]{carrillo2025longtimebehaviourbifurcationanalysis}.

We are now ready to develop in detail and apply our numerical approach to the stationary problem \eqref{eq:McKean_Vlasov_stationary}. In what follows, we fix the natural space $X=L^2(\Omega)$, the Hilbert space of (Lebesgue) measurable functions with finite $L^2$-norm. Without loss of generality, we fix $\sigma = 1$ so that $\kappa \geq 0$ is the sole bifurcation parameter to vary.

\subsection{Notations}\label{sec:notations}

We will reserve regular math letters for the analytical solution, i.e., $u(x)$ denotes the analytical solution on the continuous domain $x \in \Omega$. For simplicity, we will fix $\Omega = \mathbb{T}$ to be the torus of side length $L = \pi$. We denote by $N \gg 1$ the maximal discretisation of the domain $\Omega$, giving the fixed, maximal grid $\{ x_j \}_{j=1}^N$. This will be used when computing our reference solution for the grid refinement study. For some $m \in \mathbb{N}$ fixed, we then denote by $N_n$, $n=1, \ldots, m$, the number of grid points used. For each $1 \leq n < m$, we denote the nodes of the coarse grid as $\{ x_j ^{(n)} \}_{j=1}^{N_n}$. We will always assume uniform spacing due to the nature of operators considered here\footnote{This is not necessary in principle, but one should take care depending on the problem considered. For example, a non-uniform grid may be preferred for problems defined on an unbounded domain whose solutions are known to decay quickly in the unbounded direction (e.g., the Cucker-Smale model); on the other hand, spatial convolution on a non-uniform grid is not feasible using the FFT algorithm (e.g., the McKean-Vlasov equation).}.

Given $v \in X$, we write $\mathbf{v}^{(n)} = (v_1, \ldots, v_{N_n}) := ( v(x_j^{(n)}) )_{j=1}^{N_n}$. When an analytical solution $u(x)$ is available, we denote by $\mathbf{u}_{\textup{ref}} = ( u(x_j) )_{j=1}^N \in \mathbb{R}^N$ an analytical ``reference'' solution, which we always define on the maximal grid $\{ x_j \}_{j=1}^N$. This serves as a ground-truth reference solution to test accuracy as the grid size changes.

When we solve a given problem on a grid $\mathbf{x}^{(n)}$, we denote by $\mathbf{u}^{(n)} = ( u_j ^{(n)}  )_{j=1}^{N_n} \in \mathbb{R}^{N_n}$ the approximate solution at the nodes $\{ x_j ^{(n)} \}_{j=1}^{N_n}$. When discussing solution iterates, we write the $k^{\textup{th}}$ iterate as $\mathbf{u}^{(n),k} = ( u_{j}^{(n),k}  )_{j=1}^{N_n}$. 

Finally, for consistent comparison across mesh sizes, we upsample $\mathbf{u}^{(n)}$ with an interpolation operator $\mathcal{I}_n : \mathbb{R}^{N_n} \mapsto \mathbb{R}^N$ so that we can compare against $\mathbf{u}_{\textup{ref}}$ on the finest grid $ \{ x_j \}_{j=1}^N$. We define the upsampled solution as
$$
\widetilde{\mathbf{u}}^{(n)} := \mathcal{I}_n \mathbf{u}^{(n)} \in \mathbb{R}^N.
$$
We use the cubic interpolation implementation from the \texttt{Dierckx.jl} package \cite{dierckx1993curve, dierckxjl}.

We denote by $\mathbf{T}_n : \mathbb{R}^{N_n} \mapsto \mathbb{R}^{N_n}$ the discretisation of $\mathcal{T}$
\begin{align}\label{op:general_T_disc}
    \mathbf{T}_n ( \mathbf{v}^{(n)} ) := \left( \mathcal{T}v \left( x_j ^{(n)} \right)  \right)_{j=1}^{N_n}
\end{align}
where we numerically approximate each component of the analytical function $\mathcal{T}v(x)$. Doing this typically involves computing nonlocal component(s); for the McKean-Vlasov equation, these include the spatial convolution $W*v$ and a normalisation factor $Z(v)$. Therefore, $\mathcal{T}$ must be viewed as a nonlocal evaluation in the sense that $ \mathcal{T}v(x_j ^{(n)})$ at a single point $x_j^{(n)}$ \textit{may depend on all other elements} of $\mathbf{x}^{(n)}$.

Given a finite-dimensional operator $\mathbf{F}_n : \mathbb{R}^{N_n} \mapsto \mathbb{R}^{N_n}$ and a vector $\mathbf{v}^{(n)} \in \mathbb{R}^{N_n}$ so that $\mathbf{F}_n (\mathbf{v}^{(n)}) \in \mathbb{R}^{N_n}$, we denote the Jacobian by
$$
\mathbf{J}_{\mathbf{F}_n}[ \mathbf{v}^{(n)}] := D \mathbf{F}_n[ \mathbf{v}^{(n)}]  \in \mathbb{R}^{N_n \times N_n}
$$
In practice, we will work with the discretised operator $\mathbf{F}_n = \mathbf{T}_n - \mathbf{I}_n$ of $\mathcal{T} - I$.

As one of our goals is to verify that either the approximate or analytical Fr{\'e}chet derivative approach is viable, we distinguish between an approximate solution obtained from the analytical Fr{\'e}chet approach and an approximate solution obtained from the discrete Fr{\'e}chet approach with a superscript $\textup{AF}$ or $\textup{DF}$, i.e. $\mathbf{u}^{(n),\textup{AF}}$ or $\mathbf{u}^{(n),\textup{DF}}$. We then denote the error as a function of the chosen discretisation by
$$
\mathcal{R}_{p}^{\textup{AF}}(N_n) := \| \mathbf{u}_{\textup{ref}} - \widetilde{\mathbf{u}}^{(n),\textup{AF}} \|_p, \quad\quad \mathcal{R}_{p}^{\textup{DF}}(N_n) := \| \mathbf{u}_{\textup{ref}} - \widetilde{\mathbf{u}}^{(n),\textup{DF}} \|_p,
$$
and $\| \cdot \|_p$ is the usual $\ell^{p}$-norm in $\mathbb{R}^{N_n}$ with $p = 2$ or $\infty$.

The errors $\mathcal{R}_p^{\textup{AF}}$ and $\mathcal{R}_p^{\textup{DF}}$ assume that we have an analytical reference solution, which is not always the case. When an analytical reference solution is unavailable, we will use the numerical approximation obtained on the finest grid as the reference solution. Then, for additional grid-refinement studies of such cases, we also want to compare the use of the analytical versus discrete Fr{\'e}chet derivatives directly, in which case the error will be given by
$$
\widetilde{\mathcal{R}}_p (N_n) := \| \widetilde{\mathbf{u}}^{(n), \textup{AF}} - \widetilde{\mathbf{u}}^{(n), \textup{DF}}  \|_p.
$$
Minkowski's inequality yields the general estimate $\widetilde{\mathcal{R}}_p (N_n)  \leq \mathcal{R}_{p}^{\textup{AF}}(N_n) + \mathcal{R}_{p}^{\textup{DF}}(N_n)$.

\subsection{Discretising the McKean-Vlasov equation}
Having fixed the torus $\mathbb{T} = (-\pi/2, \pi/2]$, we now compute the nonlocal component $\kappa W*u$, where we remind readers that $\kappa \geq 0$ is an interaction strength and $W*\cdot$ denotes a periodic spatial convolution:
\begin{align}
    W*u(x) := \int_\mathbb{T} W(x-y) u(y) \textup{dy}.
\end{align}

Since we are working on the torus, we discretise the convolution $f*g$ between two periodic signals $f,g$ via a discrete convolution as described in Algorithm \ref{algorithm:periodic-convolution_intro} below.

\begin{algorithm}[!htbp]
\scriptsize
\caption{Computes the one-dimensional discrete approximation to an analytical periodic convolution for functions $f,g$, given discretisations $\mathbf{f}, \mathbf{g}$ sampled at discrete points $\mathbf{x}$, using the convolution theorem and the Fast Fourier Transform (FFT). We used the \texttt{FFTW.jl}\cite{FFTW.jl-2005} package.}
\label{algorithm:periodic-convolution_intro}
\begin{algorithmic}[1]
\Function{DiscretePeriodicConvolution}{$\mathbf{f}, \mathbf{g}, \mathbf{x}$}
    \State $\Delta x \gets x_2 - x_1$
    \Comment{Assume uniform grid.}
    \State $\mathbf{f \ast g} \gets \mathbf{0} \in \mathbb{R}^{N}$
    \State $\mathbf{h} \gets \texttt{FFT}(\mathbf{f}[1 : N - 1]) \odot \texttt{FFT}(\mathbf{g}[1 : N - 1])$
    \State $\mathbf{h} \gets \texttt{IFFT}(\mathbf{h})$
    
    \If{$\texttt{is\_symmetric\_interval}(\mathbf{x})$}
            \State $\mathbf{h} \gets \texttt{FFTShift}(\mathbf{h})$
            \Comment Perform circular shift.
    \EndIf

    \State $\mathbf{f \ast g}[1 : N - 1] \gets \Re(\mathbf{h})$
    \Comment Extract real part.
    \State $\mathbf{f \ast g}_{N} \gets \mathbf{f \ast g}_1$
    \Comment Enforce periodicity.
        
    \State \Return $\Delta x \odot \mathbf{f \ast g}$ 
    \Comment Rescale.
\EndFunction
\end{algorithmic}
\end{algorithm}

Since this can be treated using the FFT algorithm \cite{MR178586}, the speed and accuracy are a significant improvement over direct discretisation of $W * u$, especially if discretising its spatial derivatives at the level of the differential equation. In our case, we apply Algorithm \ref{algorithm:periodic-convolution_intro} to the discretised signals of the kernel $f = W$ and a candidate solution $g = u$ to obtain a discretised version of $W*u$ at each Newton iteration step.

Once $W*u$ is discretised on $\mathbf{x}$, we can perform elementwise operations to compute the discretised exponential appearing in the numerator of \eqref{eq:operator_T_McKean}. 

Lastly, the normalisation term $Z(u)$ is approximated by numerically integrating this discretised numerator over $\mathbf{x}$ using Simpson's 1/3 Rule. In some cases, such as when integrating periodic signals, trapezoidal integration may be preferred; we opt for Simpson's rule for consistency in application to each problem considered here, since not all problems are defined on a periodic domain.

We have now fully discretised the operator $\mathcal{T}$ introduced in \eqref{eq:operator_T_McKean} via $\mathbf{T}_n ( \mathbf{u}^{(n)} )$ as defined in \eqref{op:general_T_disc},where we emphasise that the elements $\left(\mathbf{T}_n ( \mathbf{u}^{(n)})\right)_j$ depend not only on the entry $u_j$, but on the full vector $(u_1, \ldots, u_{N_n})$ through the nonlocal components $W*u$ and $Z(u)$. 

Finally, we must compute the Jacobian $\mathbf{J}_\mathbf{T_n} [\mathbf{u}^{(n)}]$ and its inverse. In practice, this introduces a significant computational bottleneck as this calculation is required at each iteration. However, as discussed earlier, we have the analytical Fr{\'e}chet derivative as computed in \eqref{eq:mckean_vlasov_frechet_der}, which will reduce the impact of this bottleneck through direct evaluation of the analytical Jacobian, leaving inversion as the remaining difficulty. 

To this end, we discretise $D\mathcal{T}[u](\phi)$ by using the initial discretisation of $\mathcal{T}$ on the grid of size $N_n$, and then applying Algorithm \ref{algorithm:periodic-convolution_intro} to obtain the convolution and integrated convolution terms.

To compute the inverse, we opt for approximating $\left(\mathbf{J_{T_n}}[\mathbf{u}^{(n)}]\right)^{-1} \mathbf{T}_n(\mathbf{u}^{(n)})$ iteratively within a Krylov subspace. To do this, we employ the GMRES algorithm \cite{gmres} (from the \texttt{Krylov.jl} package \cite{krylov}) without preconditioners, which only requires being able to multiply vectors by the Jacobian. This is computed directly via:
\begin{align}\label{eq:exact_frechet_1}
    \mathbf{J_{T_n}}[\mathbf{u}^{(n)}]\boldsymbol{\phi}^{(n)} = D\mathbf{T}_n[\mathbf{u}^{(n)}](\boldsymbol{\phi}^{(n)}) \in \mathbb{R}^{N_n}.
\end{align}
This is computationally cheaper, involving only function evaluation rather than matrix multiplication, and computes exactly the evaluation of the Fr{\'e}chet derivative at $\boldsymbol{\phi}^{(n)}$. 

We are now equipped to identify the fixed points of $\mathbf{T}_n$ by applying Newton's Method to the auxiliary nonlinear map $\mathbf{F}_n(\mathbf{u}) := \mathbf{T}_n(\mathbf{u}) - \mathbf{I}_n(\mathbf{u})$. The basic algorithm is well-known, see, e.g., \cite{Kelley1995iterative}: given the map $\mathbf{F}_n : \mathbb{R}^{N_n} \to \mathbb{R}^{N_n}$ and some initial guess $\mathbf{u}^{(n),1}$ defined on the grid of size $N_n$, we iteratively perform:
\[
\forall k \geq 2, \quad \mathbf{u}^{(n),k+1} = \mathbf{u}^{(n),k} - \left(\mathbf{J_{F_n}}\left[\mathbf{u}^{(n),k}\right]\right)^{-1} \mathbf{F}_n(\mathbf{u}^{(n),k}),
\]
where $\mathbf{J_{F_n}}(\mathbf{u}^{(n),k})$ denotes the Jacobian of $\mathbf{F}_n$ at $\mathbf{u}^{(n),k}$. Algorithm \ref{algorithm:jfnk-fixed-points_intro} implements this procedure, producing an output $\mathbf{u}^{(n)}$ approximating a stationary state.

\begin{algorithm}[!htbp]
\scriptsize
\caption{Computes fixed points of a discretised operator $\mathbf{T}_n$, by finding zeros of the associated operator $\mathbf{F}_n(\mathbf{u}^{(n)}) = \mathbf{T}_n(\mathbf{u}^{(n)}) - \mathbf{I}_n\mathbf{u}^{(n)}$ with JFNK by using its the Fr{\'e}chet Derivative $D\mathbf{F}_n$. Beginning with an initial guess $\mathbf{u}^{(n)}$, iteratively refines it with $\Delta \mathbf{u}^{(n)} =  -\left(\mathbf{J_{\mathbf{T}_n}}[\mathbf{u}^{(n)}]\right)^{-1} \mathbf{T}_n(\mathbf{u}^{(n)})$, and terminates whenever the $\ell^p$-norm $\|\cdot\|_p$ of $\Delta \mathbf{u}^{(n)}$ and $\mathbf{F}_n(\mathbf{u}^{(n)})$ fall within some tolerance $\operatorname{tol} > 0$.}
\label{algorithm:jfnk-fixed-points_intro}
\begin{algorithmic}[1]
\Function{JFNKAnalyticalF{\'e}chet}{$\mathbf{T}_n : \mathbb{R}^{N_n} \to \mathbb{R}^{N_n}$, $\mathbf{u}^{(n)} \in \mathbb{R}^{N_n}$, n\_iters, tol, $p$}
    \State $\mathbf{F}(\mathbf{u}^{(n)}) := \mathbf{T}_n(\mathbf{u}^{(n)}) - \mathbf{I}_n \mathbf{u}^{(n)}$
    \State $\Delta \mathbf{u}^{(n)} \gets \mathbf{0} \in \mathbb{R}^{N_n}$
    \For{$i \in [1, \text{n\_iters}]$}
        \If{$\|\Delta \mathbf{u}^{(n)}\|_p + \|\mathbf{F}(\mathbf{u}^{(n)})\|_p \leq \operatorname{tol}$}
            \State \Return $\mathbf{u}^{(n)}$
        \EndIf
        \State $\Delta \mathbf{u}^{(n)} \gets \text{GMRES}(\mathbf{v}^{(n)} \mapsto D\mathbf{F}_n[\mathbf{u}^{(n)}](\mathbf{v}^{(n)}), -\mathbf{F}(\mathbf{u}^{(n)}))$ 
        \State $\mathbf{u}^{(n)} \gets \mathbf{u}^{(n)} + \Delta \mathbf{u}^{(n)}$
    \EndFor
    \State \Return $\mathbf{u}^{(n)}$
\EndFunction
\end{algorithmic}
\end{algorithm}

\subsection{Finite-difference approximation of the Fr{\'e}chet derivative}

In general, the multiplication in \eqref{eq:exact_frechet_1} can be approximated through finite differences. For instance, the central difference approximation at $\mathbf{v}^{(n)} \in \mathbb{R}^{N_n}$ with variation $\boldsymbol{\phi}^{(n)}$ is:
\begin{equation}
\label{eq:cd-approximation}
\boldsymbol{J_{F}}[\mathbf{v}^{(n)}]( \boldsymbol{\phi}^{(n)} ) \approx \frac{\boldsymbol{F}(\mathbf{v}^{(n)} + h \odot \boldsymbol{\phi}^{(n)}) - \boldsymbol{F}(\mathbf{v}^{(n)} - h \odot \boldsymbol{\phi}^{(n)})}{2h},
\end{equation}
where $0 < h \ll 1$ is a small parameter and $\odot$ is understood to mean elementwise multiplication. We find that the finite-difference approximation is a reasonable compromise for cases where the analytical expression is difficult or impossible to obtain (see, e.g., Figure \ref{fig:validation_McKeanVlasov} and the surrounding discussion of Section \ref{sec:results}).

\section{Results}
\label{sec:results}

We now present the results of our numerical scheme as applied to the McKean-Vlasov equation. Given the degree of freedom offered by the general input kernel $W$, we carefully choose several key examples. The first is the so-called \textit{Kuramoto model of synchronised oscillators} \cite{kuramoto1975self, kuramoto1984chemical, Carrillo2020LongTime}. This is the simplest choice, as the kernel is comprised of precisely one basis function $w_k$ defined in \eqref{eq:kernel_kuramoto} below, which produces a single stable solution branch emerging from the homogeneous state at an analytically identifiable critical threshold $\kappa^* > 0$. Increasing the complexity slightly, we then consider the same model with $W$ given by a linear combination of exactly two basis functions, $w_{k_1}$ and $w_{k_2}$. Finally, we consider some common kernels that appear in the literature, including the top-hat kernel \cite{wangsalmaniw2022} and the triangular kernel, where the latter can be obtained by convolving the former with itself, along with a third ``attractive-repulsive'' kernel made of two top-hat kernels of opposite sign. 

We remind readers that in all examples of the McKean-Vlasov equation, we solve the problem on the torus of side length $L= \pi$, subject to periodic boundary conditions as in Section \ref{sec:dev-mckean-vlasov}, discretised into 2,001 samples. All numerical experiments were performed on a MacBook Pro equipped with an Apple M2 chip and 16 GB of RAM, using Julia \cite{bezanson2017julia} version 1.11.4. Unless otherwise specified, when using Algorithm \ref{algorithm:jfnk-fixed-points_intro}, the analytical Fr{\'e}chet derivative was used, with \texttt{n\_iters} = 20,  \texttt{tol} = $1 \times 10^{-7}$ and $p = \infty$. When applicable, multithreading was enabled to accelerate computations, with up to 8 threads used concurrently using Julia’s built-in multithreading via \texttt{Threads.@threads}.

More precisely, the next subsections are devoted to verifying the following quantitative and qualitative aspects of the McKean-Vlasov model \eqref{eq:McKean_Vlasov_stationary},
the Cucker-Smale model \eqref{eq:cucker_smale_stationary} and the neural Fokker-Planck equation \eqref{eq:neural_FP_stationary}, and to explore these problems beyond the available theory:
\begin{itemize}
    \item[i)] Quantitative error estimates between a derived analytical solution and our numerical approximations as the grid size $N_n$ is increased;
    \item[ii)] Emergence of bifurcation branches, and the numerical and analytical branches obtained, agree when bifurcation analysis is available;
    \item[iii)] Produce bifurcation diagrams, demonstrating phenomena not previously observed, especially for the McKean-Vlasov equation;
    \item[iv)] A brief examination of the role of the initial iterate in the Newton scheme through a display of the basins of attraction for different initial guesses;
    \item[v)] Initial exploration of the two-dimensional McKean-Vlasov equation; identification of several inhomogeneous solutions to the Neural Fokker-Planck equation.
\end{itemize}

\subsection{The generalised Kuramoto model of synchronised oscillators} \label{sec:kuramoto_model_1}

The generalised Kuramoto model (see \cite{RevModPhys.77.137} and \cite[Section 6.1]{Carrillo2020LongTime}) considers interaction kernels of the form 
\begin{align}\label{eq:kernel_kuramoto}
    W (x) = - w_k(x),
\end{align}
where $w_k(x)$ is a basis element for $L^2_s(\mathbb{T}) \subset L^2(\mathbb{T})$ defined by $w_k(x) := \sqrt{\tfrac{2}{\pi}} \cos \left(2 k x \right)$ for any integer $k \geq 1$. Here, $L^2_s(\mathbb{T})$ denotes the (closed) subspace of $L^2(\mathbb{T})$ consisting of even functions. As $L^2_s(\mathbb{T})$ is closed in $L^2(\mathbb{T})$, it is a Hilbert space in its own right, and all Fourier coefficients of the kernel $W = - w_k(x)$ are zero except for at the single wavenumber $k$. The Kuramoto model will act as our primary validation case by following the analysis of \cite{Carrillo2020LongTime, carrillo2025longtimebehaviourbifurcationanalysis}. 

The operator of interest is then $\mathcal{T}$ defined in \eqref{eq:operator_T_McKean} with $W$ as defined in \eqref{eq:kernel_kuramoto} above. This problem has a unique homogeneous stationary state given by $u_\infty = L^{-1} = \pi^{-1}$. By \cite[Theorem 2.3 \& Proposition 2.4]{Carrillo2020LongTime}, $u^*$ is a stationary solution if and only if it is a fixed point of the map $\mathcal{T}$, or equivalently, if the condition \eqref{eq:fixed_point_preliminary} is satisfied. We now derive our analytical reference solution $\mathbf{u}_{\textup{ref}}$ as defined in Section \ref{sec:notations} against which to test our approximate solutions obtained through fixed points of the map $\mathcal{T}$.

\begin{figure}[!htbp]
\setlength{\belowcaptionskip}{0pt}
    \centering
\includegraphics[width=0.95\linewidth]{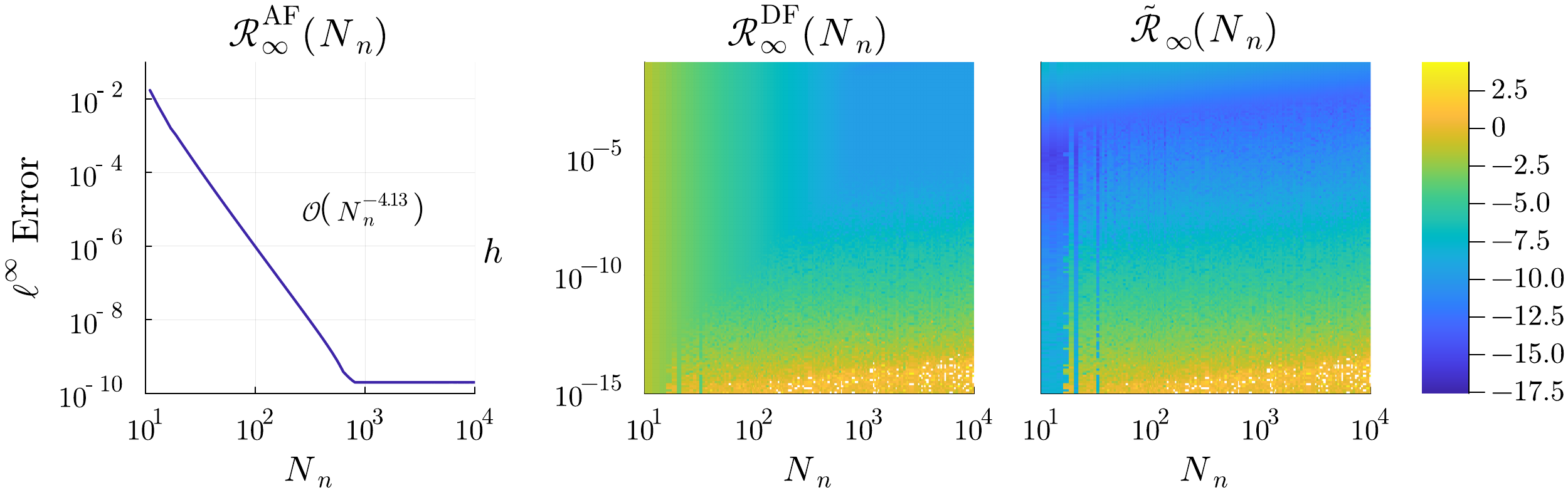}
    \caption{Error plots for the numerical solutions obtained to the McKean-Vlasov equation \eqref{eq:McKean_Vlasov_stationary} with kernel $W = - w_k(x)$ as developed in Section \ref{sec:kuramoto_model_1} at the fixed value $\kappa = 3$. For the left panel, we use 2,017 logarithmically spaced (odd) values of $N_n \in [11, 10001]$. Note that due to the logarithmic spacing, the first values of $N_n$ will be consecutive (i.e $N_1 = 11, N_2 = 13, N_3 = 15$). For the centre and right panels, we employ 160 logarithmically spaced (odd) values of $N_n \in [11, 10001]$ and $h \in [10^{-15}, 10^{-1}]$. As an initial guess for Algorithm \ref{algorithm:jfnk-fixed-points_intro}, $f(x) = 1/\pi + \cos(2x)$ was used. Any white points in the heatmaps correspond to instances in which Algorithm \ref{algorithm:jfnk-fixed-points_intro} failed to converge when using the central-differences approximation of the Fr{\'e}chet derivative.}
    \label{fig:validation_McKeanVlasov}
\end{figure}

\begin{figure}[!htbp]
\setlength{\belowcaptionskip}{0pt}
    \centering
\includegraphics[width=0.8\linewidth]{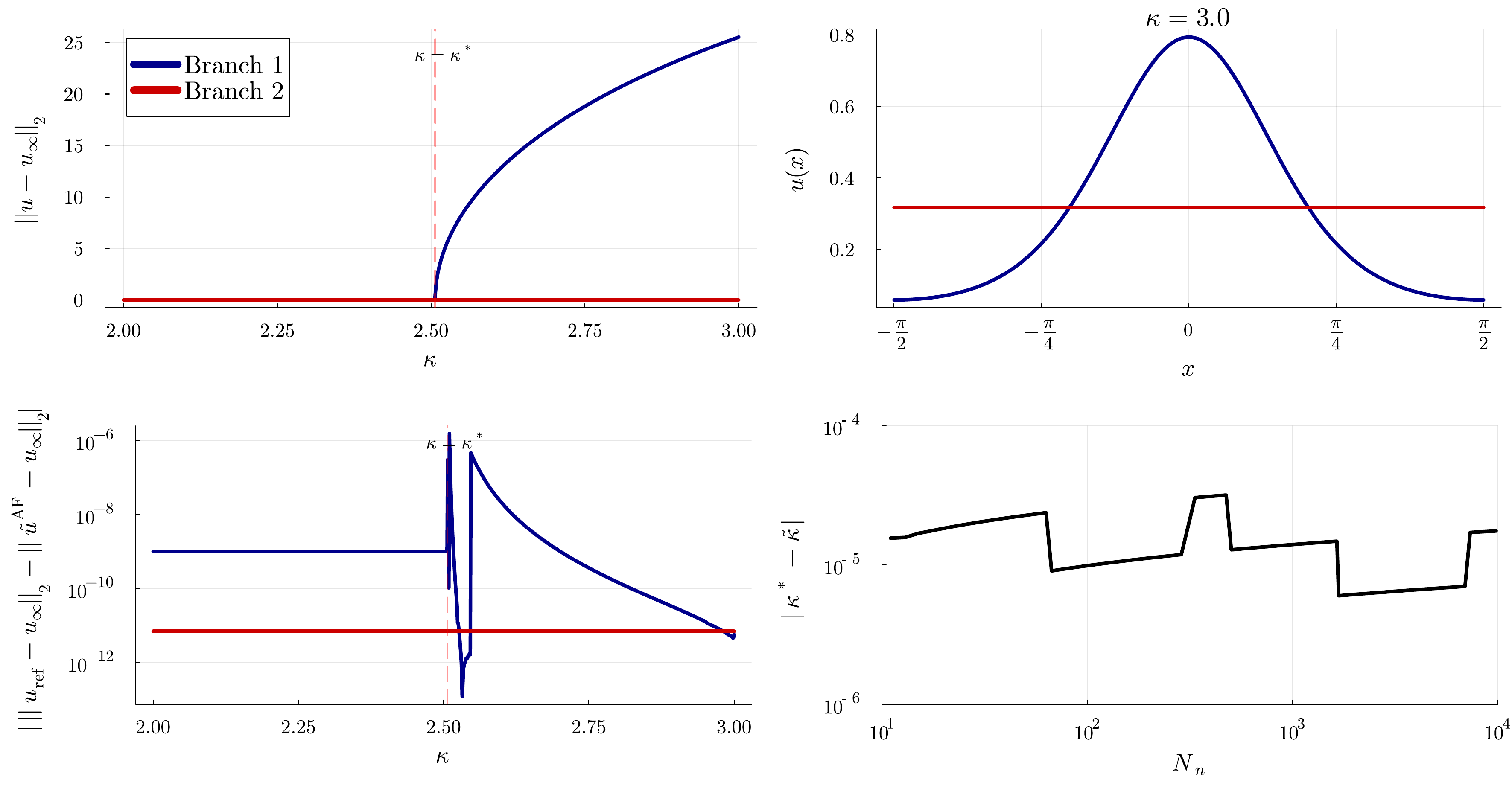}
    \caption{Top-left panel: bifurcation diagram produced using our algorithm. Top-right panel: the two solution profiles (homogeneous and single nontrivial solution) at $\kappa = 3$. Bottom-left panel: the error between the bifurcation diagrams obtained using the analytical reference solution and our numerical approximation. Bottom-right panel: the error between the critical analytical threshold $\kappa^*$ and that predicted by our numerical approximations, across 160 logarithmically spaced (odd) values of $N_n \in [11, 10001]$. 
    }
    \label{fig:validation_bifurcation_McKeanVlasov}
\end{figure}

\noindent\textbf{Derivation of the reference solution.} We use the approach of \cite[Sec. 6.1]{Carrillo2020LongTime}. To find the nontrivial analytical solution, we first identify the fixed point $a = a_{\kappa}$ of
\begin{align}
        M(a, \kappa) = \sqrt{\frac{2}{\pi}} \kappa \frac{I_1(a)}{I_0(a)},
    \end{align}
where $I_k, k \in \mathbb{N}$ is the modified Bessel function of the first kind of order $k$ \cite{MR1349110}. The analytical solution for the Kuramoto model is then given by    
\begin{align}
        \label{eq:analytical-kuramoto}
        u_{\textup{Kuramoto}}(x) = \frac{1}{\pi}\frac{\exp\left(a_\kappa \cos\left(2 k x\right)\right)}{I_0(a_\kappa)}.
    \end{align}
Using Julia's \texttt{SpecialFunctions.jl} and Newton's method, we are able to identify the fixed point $a_\kappa$ to machine precision, that is, so that $| M(\tilde a_\kappa, \kappa) - \tilde a_\kappa | = 0$ when computed inside Julia, where $\tilde a_\kappa$ is our numerical approximation of $a_\kappa$. Our analytical reference solution is then obtained by discretising \eqref{eq:analytical-kuramoto} on the finest grid:
\begin{align}\label{eq:u_ref_kuramoto}
    \mathbf{u}_{\textup{ref}} := ( u_{\textup{Kuramoto}}(x_j) )_{j=1}^N
\end{align}
Using this analytical reference solution, we consider two circumstances as validation of our approach: first, the quantitative error estimates at a fixed value of $\kappa$ (Figure \ref{fig:validation_McKeanVlasov}); second, further quantitative and qualitative aspects across a range of $\kappa$ values (Figure \ref{fig:validation_bifurcation_McKeanVlasov}). We explain these two cases in more detail now.

\noindent\textbf{Error estimates.} From the reference solution $\mathbf{u}_{\textup{ref}}$ obtained in \eqref{eq:u_ref_kuramoto}, we can directly compare with our numerically obtained solutions, namely $\widetilde{\mathbf{u}}^{(n),\textup{AF}}$ (using analytical Fr{\'e}chet) or $\widetilde{\mathbf{u}}^{(n),\textup{DF}}$ (using discrete Fr{\'e}chet via central difference), depending on the number of spatial points $N_n$, after (cubic) interpolation of the numerical solutions from a grid of size $N_n$ onto the finest grid of size $N$.

In Figure \ref{fig:validation_McKeanVlasov}, we display the error estimates at the fixed value $\kappa = 3$ (note that from a linear stability analysis, we know the critical threshold is $\kappa^* = \sqrt{2 \pi} \approx 2.5$, and so an inhomogeneous solution is known to exist at $\kappa = 3 > \kappa^*$). In the left-most panel of Figure \ref{fig:validation_McKeanVlasov}, we then display the error $\mathcal{R}_\infty ^{\textup{AF}}( N_n)$, as a function of the number of discretisation points $N_n$. We find that the error decays like $\mathcal{O}(N_n^{-4})$. 

When using a finite-difference approximation of the Fr{\'e}chet derivative, we also have the parameter $h$ to consider; therefore, in the centre panel of Figure \ref{fig:validation_McKeanVlasov}, we display the error $\mathcal{R}_\infty ^{\textup{DF}}( N_n)$ as a heat map, depending on the number of discretisation points $N_n$ and the size of the small parameter $h$ used in the central-difference appoximation of the Fr{\'e}chet derivative as defined in \eqref{eq:cd-approximation}. For moderate values of $h$ (approximately $10^{-1}$ to $10^{-7}$), horizontal cross-sections recover the same error decay rate as displayed in the left-most panel, i.e., $\mathcal{R}_\infty ^{\textup{DF}}( N_n)$ also decays like $\mathcal{O}(N_n^{-4})$; when $h$ is chosen too small ($h \ll 10^{-7}$), we observe diminishing returns as $N_n$ increases. 

Finally, in the right-most panel of Figure \ref{fig:validation_McKeanVlasov} we display the error $\mathcal{R}_\infty( N_n)$ directly between the analytical and discrete implementation of the Fr{\'e}chet derivative. The error follows a similar pattern to the left and centre panels described above. For small values of $N_n$ (i.e., $N_n \sim 10^{1}$), we observe very small error; we carefully note that this does not suggest the solution obtain is a good one, but rather that that two solutions obtained are close to each other, regardless of the size of $\norm{\mathcal{T}u - u}$. 

Altogether, we have verified that: the scheme can recover the analytical solution to great precision;
the error converges approximately on the order $N_n^{-4}$;
and that the finite-difference approximation of the Fr{\'e}chet derivative is effective in the absence of an analytical version; using thediscrete or analytical Fr{\'e}chet derivative recovers the true solution with great precision, provided that the parameter $h$ is chosen appropriately.

Based on our explorations, these findings do not appear to depend on this specific instance of the McKean-Vlasov equation; the same convergence properties hold across all kernels tested, with the caveat that when it cannot be derived directly, the approximate solution obtained on the finest grid acts as the reference solution.

\noindent\textbf{Further qualitative validation.} With verification at a fixed value of $\kappa$, we then ran the solver sweeping across $\kappa \in [2,3]$ to recover the bifurcation diagram of \cite[Figure 2]{Carrillo2020LongTime}, which we display in the top-left panel\footnote{Note carefully that we display the norm difference $\norm{u - u_\infty}_2$ on the vertical axis; \cite[Figure 2]{Carrillo2020LongTime} uses a different, but equivalent, metric on the vertical axis. The qualitative features of the diagram remain the same, but the true numerical values will differ.} of Figure \ref{fig:validation_bifurcation_McKeanVlasov}. We now explain how we produce the top-left panel of Figure \ref{fig:validation_bifurcation_McKeanVlasov}, and how we produce our subsequent bifurcation diagrams. We first run Algorithm \ref{algorithm:jfnk-fixed-points_intro} to identify an inhomogeneous solution at $\kappa = 3 > \kappa^*$; this solution is found in the top-right panel of Figure \ref{fig:validation_bifurcation_McKeanVlasov}. We numerically validate this solution using the steady state condition \eqref{eq:fixed_point_preliminary} obtaining an error consistent with the tolerance of the method. The error in satisfying condition \eqref{eq:fixed_point_preliminary} is measured as the maximum of the absolute error of the left-hand side of \eqref{eq:fixed_point_preliminary} with respect to its average.

Once a nontrivial solution is identified, this becomes the initial guess at $\kappa = 3 - \delta$ for a small step $\delta > 0$. We continue this procedure to trace the entire branch down to $\kappa = \kappa^*$; once the homogeneous branch is reached at $\kappa = \kappa^*$, subsequent iterations for $\kappa < \kappa^*$ will yield the homogeneous solution only. In this instance, we know analytically that there is exactly one branch to identify; therefore, once an inhomogeneous state has been identified, we need not continue searching for others. This is in contrast to the general case, where initial guesses may need to be chosen carefully, and sweeping in the $\kappa$ domain may need to be repeated in both directions to trace all branches (see Figure \ref{fig:two_mode_bifurcation} and the discussion of Section \ref{sec:init_iterate}). We also numerically validate the steady state condition \eqref{eq:fixed_point_preliminary} over the entire bifurcation branch(es) obtaining again an error consistent with the tolerance of the method.

In the bottom-left panel of Figure \ref{fig:validation_bifurcation_McKeanVlasov}, we display the difference between the bifurcation curves produced by the reference solution \eqref{eq:u_ref_kuramoto} and the numerical solution obtained via the analytical Fr{\'e}chet derivative as a function of $\kappa$. We observe that their difference does not exceed $\sim 10^{-6}$ across all values of $\kappa \in [2,3]$. This can be viewed as similar verification found in the left panel of Figure \ref{fig:validation_McKeanVlasov}, but now across a range of $\kappa$.

In the bottom-right panel of Figure \ref{fig:validation_bifurcation_McKeanVlasov}, we display the difference between the analytical and numerical values of the critical threshold $\kappa^*$; this is what is found visually in the top-left panel, but we also verify this quantitatively as follows. First, we identified a $\kappa$ range such that the left-most bound returns a homogeneous solution (e.g., $\kappa = 2$), while the right-most bound returns an inhomogeneous solution (e.g., $\kappa = 3$). We then apply a standard bisection method over a range of $\kappa$ values such that this property is retained; once the left and right bounds are within some predefined tolerance, we terminate the bisection procedure. In the bottom-right panel of Figure \ref{fig:validation_bifurcation_McKeanVlasov}, we observe that this approximation of the critical threshold is accurate up to $\sim 10^{-5}$ across all values $N_n$, even for relatively small discretisation sizes $N_n \sim 10^1$. 

Altogether, we have verified that: the numerical scheme produces a single solution branch emerging close to this critical threshold $\kappa = \kappa^*$ as a supercritical bifurcation (see \cite[Section 6.1]{carrillo2025longtimebehaviourbifurcationanalysis} and \cite[Theorem 1.2]{carrillo2025longtimebehaviourbifurcationanalysis}); the emergent solution is at frequency $k=1$, corresponding to the single-mode kernel; and the solver can accurately recover the critical threshold $\kappa^*$, even for moderate discretisation sizes $N_n$.

\subsection{The bimodal kernel}\label{sec:two_mode_mckean_vlasov}

We now increase the complexity slightly from the Kuramoto model to the following two-mode modification:
\begin{align}\label{eq:two-mode-kernel}
    W = -a_1 w_{k_1} - a_2 w_{k_2},
\end{align}
where $a_i > 0$, $i=1,2$, are fixed constants to be chosen appropriately. By analysing the spectrum of the linearised operator, we can identify two points of criticality in the $\kappa$ domain:
\begin{align}
    \kappa_i := \frac{\sqrt{2\pi}}{a_i}, \quad i = 1,2.
\end{align}
According to \cite{Carrillo2020LongTime, carrillo2025longtimebehaviourbifurcationanalysis}, a nontrivial solution of the form
\begin{align}
\label{eq:nontrivial-two-mode-solution}
u_i(x) = \frac{1}{\pi } + s w_{k_i} + o(s^2), \quad \as{s} \ll 1,
\end{align}
emerges near $\kappa = \kappa^i$, where $w_{k_i}$ is an orthonormal basis element as defined at the beginning of Section \ref{sec:kuramoto_model_1}. Therefore, we anticipate the following:
\begin{itemize}
    \item[i)] there are (at least) two nontrivial stationary states $u_i$ to identify numerically;
    \item[ii)] $u_i$ should emerge near their respective values $\kappa = \kappa_i$;
    \item[iii)] near the trivial solution $u_\infty = 1/\pi$, $u_i$ should emerge with frequency $k_i$.
\end{itemize}

\noindent Furthermore, based on the recent results of \cite[Theorem 1.2]{carrillo2025longtimebehaviourbifurcationanalysis}, we have a precise analytical bifurcation structure that we use to verify our numerical observations further. Different from the single-mode case, a subcritical bifurcation is now possible: when $k_2 = 2 k_1$, $k_1 \geq 1$, there is a window in which a subcritical bifurcation occurs.

\begin{figure}[!htbp]
\setlength{\belowcaptionskip}{0pt}
  \centering
  \begin{subfigure}[b]{\textwidth}
    \centering
    \includegraphics[width=0.8\linewidth]{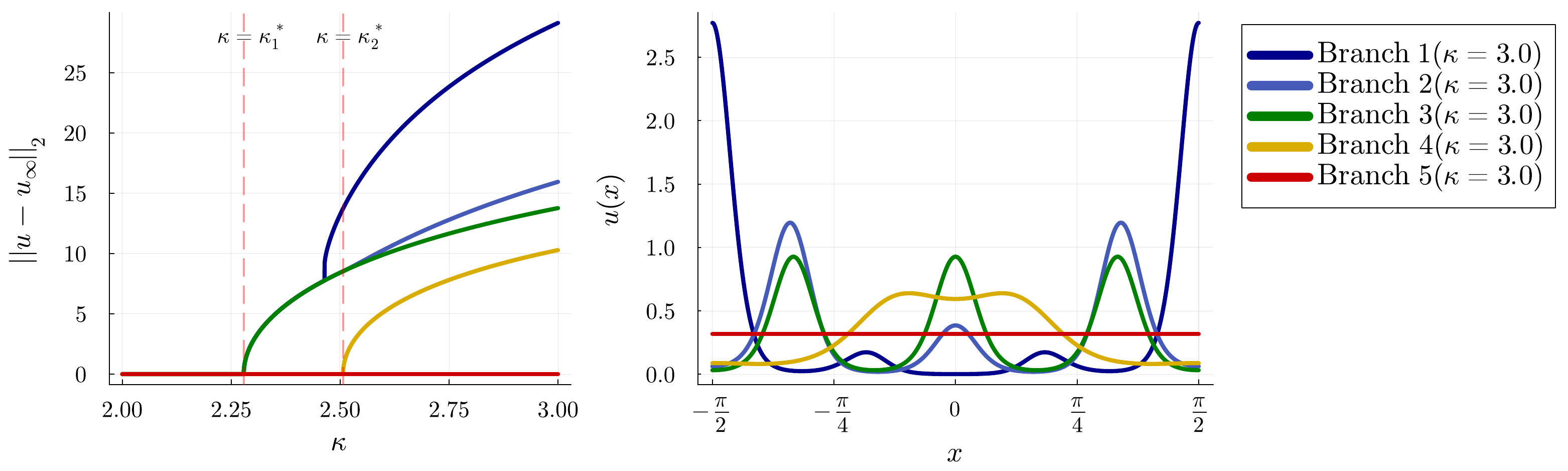}
    \caption{Supercritical bifurcation, with $a_1 = 1.0, a_2 = 1.1, k_1 = 1, k_2 = 3$, computed in $6.267 \operatorname{s} \pm 368.276 \operatorname{ms}$.}
    \label{fig:supercritical_bifurcation}
  \end{subfigure}
  \begin{subfigure}[b]{\textwidth}
    \centering
    \includegraphics[width=0.8\linewidth]{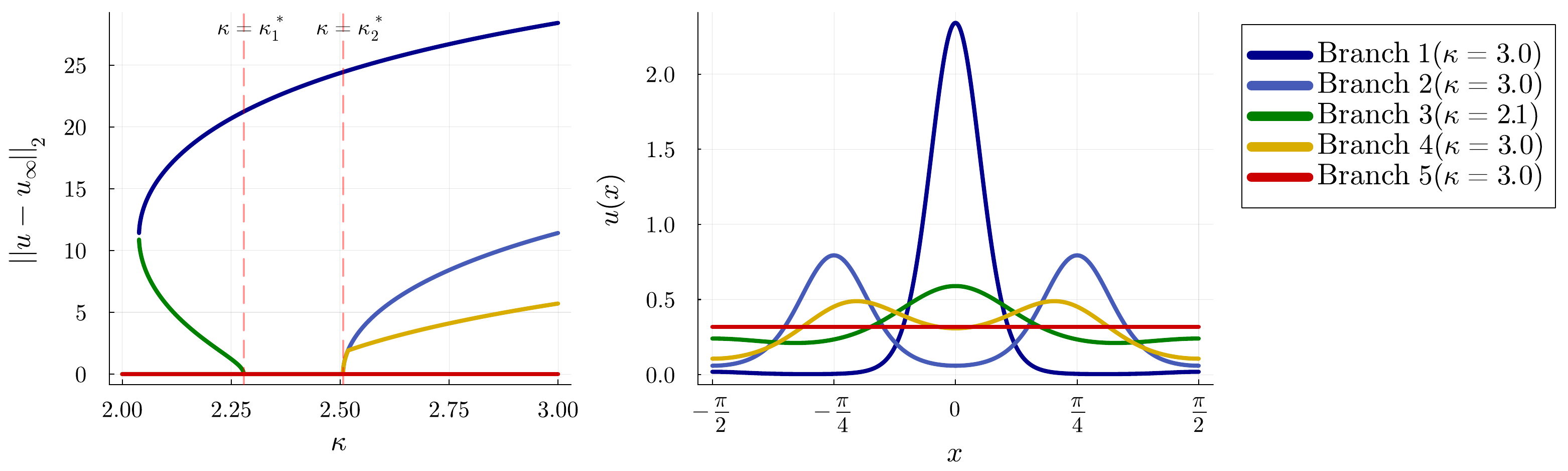}
    \caption{Subcritical bifurcation, with $a_1 = 1.1, a_2 = 1, k_1 = 1, k_2 = 2$, computed in $7.114 \operatorname{s} \pm 218.468 \operatorname{ms}$.}
    \label{fig:subcritical_bifurcation}
  \end{subfigure}
  \begin{subfigure}[b]{\textwidth}
    \centering
    \includegraphics[width=0.8\linewidth]{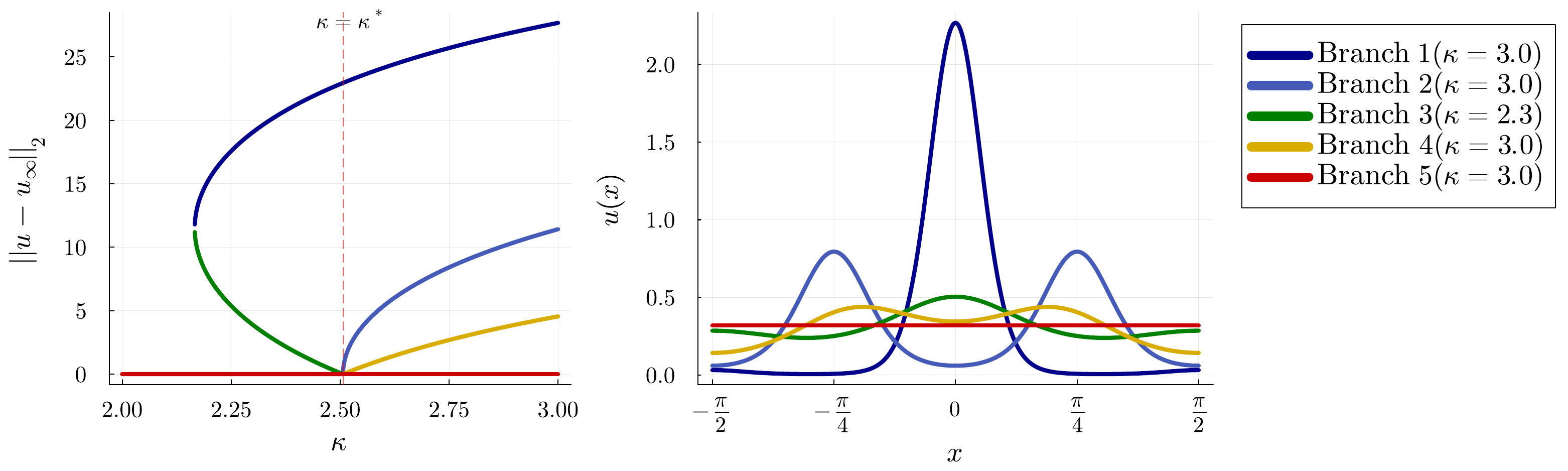}
    \caption{Bifurcation, with $a_1 = a_2 = 1, k_1 = 1, k_2 = 2$, computed in $7.425 \operatorname{s} \pm 336.020 \operatorname{ms}$}
    \label{fig:equal_coefficients_bifurcation}
  \end{subfigure}

  \caption{Bifurcation curves (left panels) and solution profiles (right panels) for the two-mode kernel defined in \eqref{eq:two-mode-kernel} for problem \eqref{eq:McKean_Vlasov_stationary}. We consider $2001$ evenly spaced values of $\kappa \in [2,3]$, and $2001$ evenly spaced values of the domain $\mathbb{T} = (-\pi/2, \pi/2]$. We report the average and standard deviation of the run times, across 10 runs.}
  \label{fig:two_mode_bifurcation}
\end{figure}

To compute the bifurcation curves shown in Figure \ref{fig:two_mode_bifurcation}, we developed a more general and robust approach than the one used for Figure \ref{fig:validation_bifurcation_McKeanVlasov}. For general kernels or PDEs, the stationary solutions are often not known in closed form, making it difficult to choose good initial guesses for Algorithm \ref{algorithm:jfnk-fixed-points_intro}. To address this, we first perform a sweep over initial guesses to identify distinct stationary states at $\kappa = 3 > \kappa^*$. For example, in this instance, building from \eqref{eq:nontrivial-two-mode-solution} we generate 91 candidate initial guesses of the form
\[
u(x;z) = \frac{1}{\pi} + 0.5 \cos((1 + 0.1z)x), \qquad z \in [0,90],
\]
along with the homogeneous state $u_\infty = 1/\pi$.
Each $z$ value then produces an input to Algorithm \ref{algorithm:jfnk-fixed-points_intro}, and the resulting stationary states are filtered to remove duplicates (up to periodic shifts). The remaining solutions form the initial set of branches used for continuation in $\kappa$. We again numerically validate the steady state condition \eqref{eq:fixed_point_preliminary} over the entire bifurcation branch(es), obtaining again an error consistent with the tolerance of the method, as described in Section \ref{sec:kuramoto_model_1}.

Starting from the states found at $\kappa = 3$, we decrease $\kappa$ incrementally from $3$ to $2 < \kappa^*$, using the approximate solution obtained at the end of each implementation of the iteration scheme as the initial guess for the next one. This continuation yields the bifurcation diagram shown in Figure \ref{fig:supercritical_bifurcation}.
For Figure \ref{fig:subcritical_bifurcation} and Figure \ref{fig:equal_coefficients_bifurcation}, the same procedure is repeated, but we now use the full set of initial guesses at several values of $\kappa$ (not just the largest one) so that we may detect additional bifurcating branches. To improve efficiency, we parallelise across eight threads, with each thread responsible for performing the continuation procedure on a single branch.

We note that in Figure \ref{fig:two_mode_bifurcation}, and in all such subsequent figures, branches are ordered based on the largest value of $\|u - u_\infty\|_\infty$, with Branch 1 having the largest such value. When plotting solution profiles for the different branches, we selected the solutions computed at the largest value of $\kappa$ available; when this is not possible, such as with Branch 3 in Figure \ref{fig:subcritical_bifurcation} or Figure \ref{fig:equal_coefficients_bifurcation}, a solution computed at an intermediate value of $\kappa$ is chosen which showcases the properties of the solution profile well.

In Figure \ref{fig:two_mode_bifurcation}, we observe precisely that behaviour expected as a result of the bifurcation analysis results. In Figure \ref{fig:supercritical_bifurcation}, we find exactly two supercritical branches emerging very close to the expected bifurcation values $\kappa_i^*$, as predicted by the bifurcation analysis of \cite{carrillo2025longtimebehaviourbifurcationanalysis}. Interestingly, we also identify two additional branches emerging from the first bifurcation branch that are not explained by the existing theory.

In Figure \ref{fig:subcritical_bifurcation}, we observe a subcritical bifurcation (the first branch at $\kappa = \kappa_1^*$), and then a supercritical branch (the second branch at $\kappa = \kappa_2^*$). This is again consistent with the bifurcation analysis of \cite{carrillo2025longtimebehaviourbifurcationanalysis}, especially for the case where a subcritical bifurcation is possible for a relatively narrow range of parameter values.

Finally, in Figure \ref{fig:equal_coefficients_bifurcation}, we display the result when the coefficients $a_i$ are equal. This is a particular case that the analytical theory cannot easily handle, as it corresponds with a two-dimensional kernel of the linearised problem, and the theory of Crandall-Rabinowitz (sometimes referred to as \textit{bifurcation from a simple eigenvalue}) does not immediately apply. Our numerical results suggest that a single supercritical branch still appears at the critical threshold $\kappa_1^* = \kappa_2^* = \kappa$, alongside two additional branches that appear to be transcritical from visual inspection. This leaves open an interesting direction of study for higher-dimensional kernels (in the PDE operator sense).

\subsection{The general McKean-Vlasov Equation}\label{sec:general_MVE}
In this subsection, we apply the algorithm to general cases. 
Under some technical criteria, it is known that every wavenumber $k \in \mathbb{N}$ such that $\widetilde{W}(k) < 0$ leads to a bifurcation from the trivial solution $1/\pi$ near $\kappa = \kappa_k^*$ defined by
\begin{align}\label{eq:critical_kappa_general}
    \kappa_k^* := -\frac{\sqrt{2 \pi}}{\widetilde{W}(k)} > 0,
\end{align}
and the solution branch emerges with frequency $k$. This means that for a general kernel, we can analytically determine a minimal number of stationary states we expect to find, along with their associated frequency near the point of bifurcation. As suggested by the increase in complexity from the single-mode to the two-mode cases explored earlier, general kernels exhibit many nontrivial solutions that emerge as secondary bifurcations from those predicted by the local theory. Again, using the bifurcation analysis of \cite{carrillo2025longtimebehaviourbifurcationanalysis}, we can predict when these local bifurcations are sub- or supercritical, and all instances we explored (including those not presented here) are consistent with this fact. Rather than exploring further verification, we now display some of the interesting solution behaviour that this approach allowed us to discover.
To this end, we introduce the (attractive) top-hat kernel defined as
\begin{align}\label{kernel:tophat}
W_{\textup{tophat}} (x) := \begin{cases}
    -\frac{1}{2R}, \qquad \as{x} \leq R; \cr 
    0, \qquad \text{otherwise}.
\end{cases}
\end{align}
and the (attractive) triangle-hat kernel defined as
\begin{align}\label{kernel:triangle}
W_{\textup{triangle}} (x) := \begin{cases}
    -\frac{1}{2R} \left(1 - \frac{|x|}{R}\right), \qquad \as{x} \leq R; \cr 
    0, \qquad \text{otherwise}.
\end{cases}
\end{align}

\begin{figure}[!htbp]
\setlength{\belowcaptionskip}{0pt}
  \centering

  \begin{subfigure}[b]{\textwidth}
    \centering
    \includegraphics[width=0.8\linewidth]{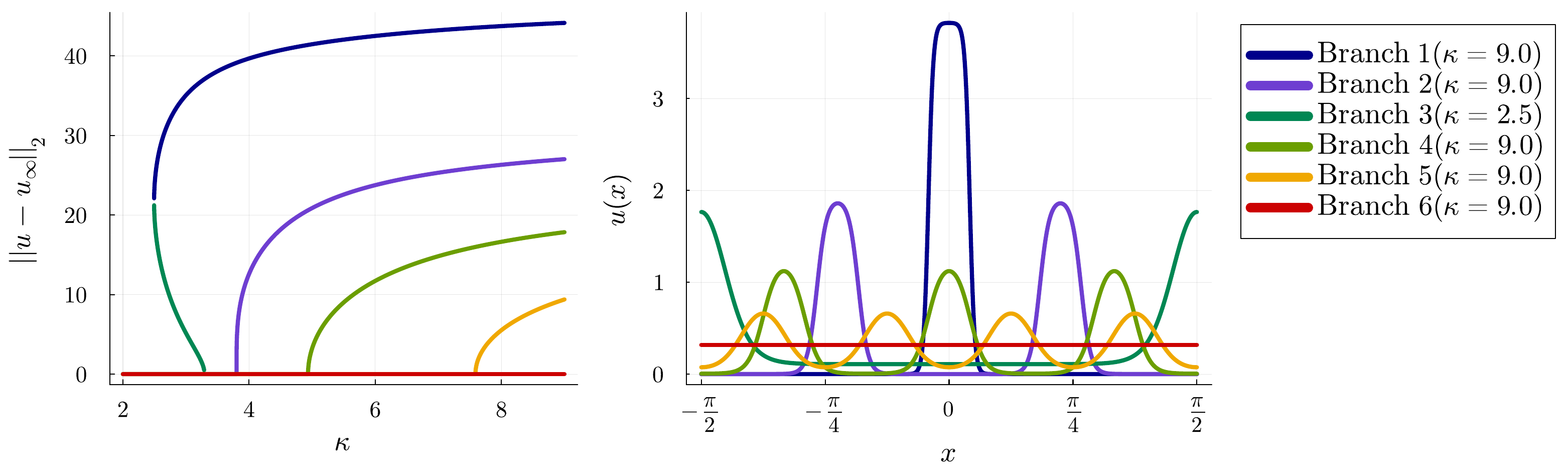}
    \caption{Bifurcation diagram using the tophat kernel \eqref{kernel:tophat} with $R = \frac{\pi}{12}$, and $\kappa \in [2,9]$.}
    \label{fig:top_hat_bifurcations}
  \end{subfigure}

\begin{subfigure}[b]{\textwidth}
    \centering\includegraphics[width=0.8\linewidth]{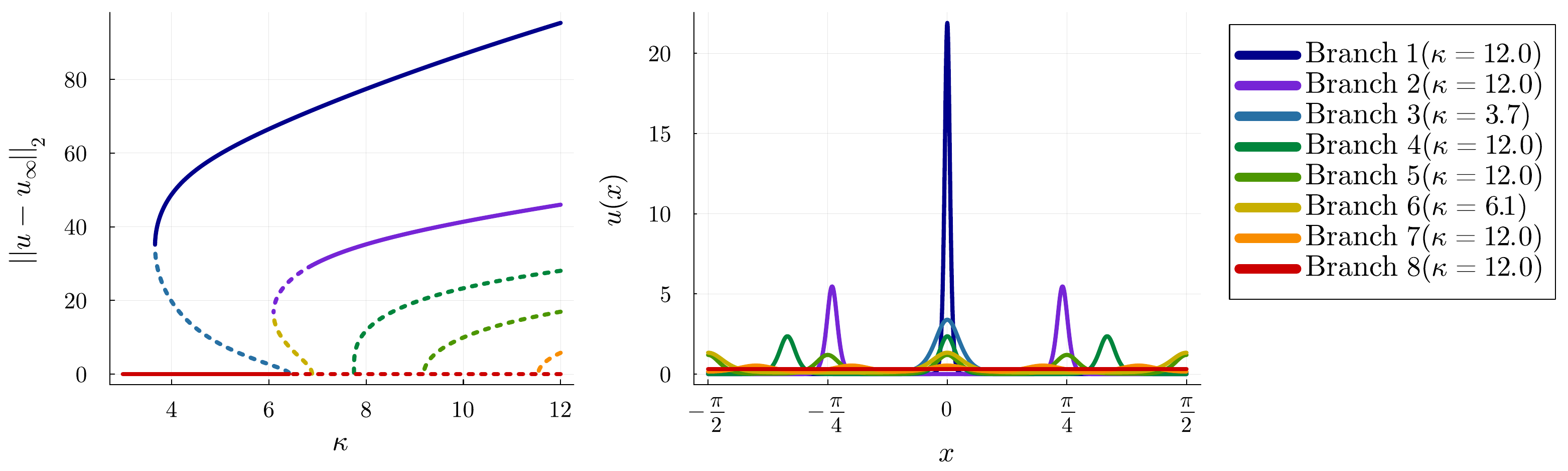}
    \caption{Bifurcation diagram using the triangular-hat kernel \eqref{kernel:triangle} with $R = \frac{\pi}{12}$, and $\kappa \in [3,12]$. For this figure only, solid lines indicate local stability of the branch, while dashed lines indicate instability (see Section \ref{sec:branch_stability}).}
    \label{fig:triangle_bifurcations}
\end{subfigure}

\begin{subfigure}[b]{\textwidth}
\centering\includegraphics[width=0.8\linewidth]{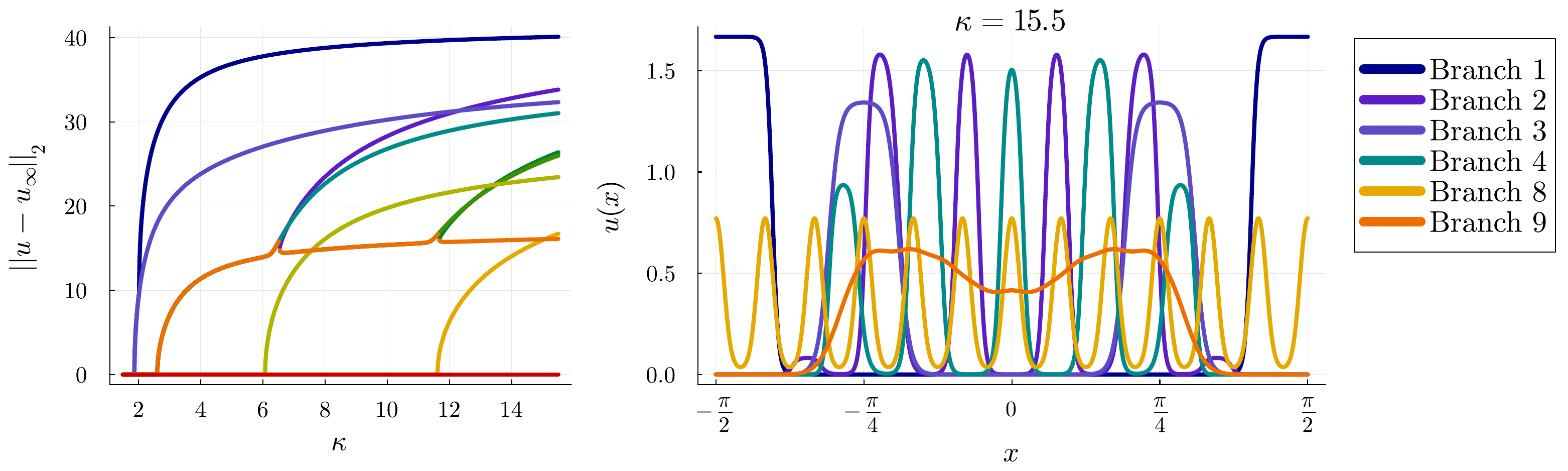}
    \caption{Bifurcation diagram using the attractive-repulsive tophat kernel \eqref{kernel:att-rep-tophat} with $R = 0.6$, and $\kappa \in [1.5,15.5]$. We display only 6 of the 10 profiles in the right panel to improve readability.}
    \label{fig:att_rep_bifurcations}
\end{subfigure}
\caption{Bifurcation curves obtained for three different instances of problem \eqref{eq:McKean_Vlasov_stationary}. We consider $2001$ evenly spaced values of the domain $\mathbb{T} = (-\pi/2, \pi/2]$ for (A), (B), and $4001$ samples for (C).}\label{fig:wacky_bifurcations}
\end{figure}

\noindent These two kernels are exemplary in a few ways. First, the tophat kernel has a relatively simple spectrum and can be computed analytically:
$$
\widetilde W_{\textup{tophat}} (k) \sim \sinc (R k).
$$
The Fourier modes then oscillate so that both positive and negative Fourier modes of $W_{\textup{tophat}}$ exist; the negative modes will produce a bifurcation branch according to \eqref{eq:critical_kappa_general}, while the positive modes do not produce a bifurcation (from the trivial branch). In Figure \ref{fig:top_hat_bifurcations}, we display the first $4$ branches for the tophat kernel identifiable via a local bifurcation analysis when $R = \pi/12$.

One can then view the modes of $W_{\textup{triangle}}$ through the convolution theorem: up to scaling, there holds $W_{\textup{tophat}} * W_{\textup{tophat}} \sim W_{\textup{triangle}}$. Hence,
$$
\widetilde{W}_{\textup{triangle}} (k) \sim - \as{\widetilde W_{\textup{tophat}} (k)}^2 \sim -\as{\sinc (R k)}^2.
$$
Consequently, except for some exceptional cases, the triangle kernel yields a bifurcation branch at \textit{every} wavenumber $k\geq1$. Moreover, since $\widetilde{W}_{\textup{triangle}} (k) \sim (kR)^{-2}$, through formula \eqref{eq:critical_kappa_general} we expect the bifurcation values $\kappa_i^*$ for the triangle kernel to be much larger than those for the tophat kernel, as $\widetilde W_{\textup{tophat}} (k) \sim (kR)^{-1}$. This is what we observe in Figure \ref{fig:triangle_bifurcations}: every wavenumber produces a bifurcation branch, and we identify the first $5$ of them. Compared to the bifurcation points found in Figure \ref{fig:top_hat_bifurcations}, the bifurcation values in Figure \ref{fig:triangle_bifurcations} are moderately larger.

These two observations motivate the following ``attractive-repulsive'' tophat kernel:
\begin{align}\label{kernel:att-rep-tophat}
W_{\textup{att-rep-tophat}} (x) := \begin{cases}
    -\frac{1}{2R}, \qquad \as{x} \leq R; \cr 
    \frac{1}{R}, \qquad R < \as{x} \leq 2R; \cr
    0, \qquad \text{otherwise}.
\end{cases}
\end{align}
This pairs two behaviours that, together, are very different from the tophat or triangle kernel: it features attractive-repulsive interactions and has the low-regularity properties of the tophat kernel. In Figure \ref{fig:att_rep_bifurcations}, we display the branches identified when $R = 0.6$. We observe that, unlike the simpler tophat and triangle kernels, there are secondary bifurcations, as found in the two-mode kernel case. We also observe several points of intersection among the identified curves. We remark that at these points of intersection, the solution profiles appear to be distinct. For example, consider the intersection point between Branch 8 and Branch 9 in Figure \ref{fig:att_rep_bifurcations}: the corresponding profiles in the right panel differ substantially, even though they are close to the $\kappa$ value at which the two curves intersect in the bifurcation diagram.

\subsection{Branch stability}\label{sec:branch_stability}
Solving directly for stationary profiles allows one to identify any solution profile, independent of the solution's stability properties. Consequently, we have no stability information beyond the local analysis of \cite{Carrillo2020LongTime, carrillo2025longtimebehaviourbifurcationanalysis}. More precisely, the point of critical stability $\kappa^*$ for the homogeneous solution (i.e., $\kappa^*>0$ such that $u_\infty = L^{-1}$ is locally stable for $\kappa < \kappa^*$, and is unstable for $\kappa > \kappa^*$) is given by the smallest positive $\kappa_k^*$ defined in \eqref{eq:critical_kappa_general}. If the branch emerging from the point of critical stability is supercritical, its solution is locally stable in $\kappa \in (\kappa^*, \kappa^*+\delta)$ for some $\delta$ small enough; if it is subcritical, the branch is unstable near $\kappa^*$. 

To explore stability properties further, we employ the numerical scheme developed in \cite{MR3372289, MR4605931}. This method preserves several desirable properties of the continuous problem (positivity, boundedness, mass, and energy decay) and is suitable for assessing the stability of the profiles we generate by introducing a small, mass-preserving perturbation to the stationary profile, and using the perturbed solution as the initial data. This is what is displayed in Figure \ref{fig:triangle_bifurcations}: dashed lines correspond to numerical instability; solid lines correspond to numerical local stability. First, we note some analytical properties that are guaranteed. From the bifurcation analysis of \cite{Carrillo2020LongTime}, we know that the homogeneous state (red line, Branch 8) is stable (solid) until $\kappa \approx 6.5$; beyond this, it is unstable (dashed). Similarly, from the bifurcation analysis of \cite{carrillo2025longtimebehaviourbifurcationanalysis}, the first emerging branch (light-blue, Branch 3) is subcritical and therefore unstable nearby.

To the best of our knowledge, analytical descriptions of local stability do not exist for all other cases. We sample $\kappa$ in intervals of length $0.25$, giving between $6$ and $24$ points per branch, and systematically add varying levels of noise to the stationary profiles identified by Algorithm \ref{algorithm:jfnk-fixed-points_intro}: large noise (relative $L^2$-error of $10^{-2}$), moderate noise (relative $L^2$-error between $10^{-3}$ and $10^{-4}$) and small noise (relative $L^2$-error of $10^{-5}$). In all simulations, we maintain the spatial discretisation of the Newton scheme ($N = 2001$ grid points) and use an implicit time-stepping algorithm (unconditionally stable) with a time step of $\rm{d}t = 0.1$. We run each simulation until a final time $T=14$, which appeared sufficiently long to capture stability/instability properties.

Across all noise levels, Branch 1 is consistently numerically stable, even near the turning point. While we do not aim to make claims about global stability, we note that most unstable solutions identified converge to the solution on Branch 1.

Similarly, across all noise levels, branches 3 and 6 (subcritical) are consistently unstable: even for small noise, the dynamics quickly move away from them. This is consistent with the local analysis corresponding to Branch 3. Branches 4, 5, and 7 (supercritical) are also reliably unstable, quickly moving off the branch even for small noise levels. 

Branch 2 was more subtle. For large values of $\kappa$ ($\kappa \approx 8.5$ or larger), it appears locally stable across all noise levels. For smaller values of $\kappa$ ($\kappa \approx 8$ or smaller), the branch appears locally stable up to $T=14$ only when adding small enough noise. The reduction of noise necessary to observe this behaviour appears to be related to how close one is to the turning point (where Branch 2 meets Branch 6 near $\kappa \approx 6$). Under even smaller noise (relative $L^2$-error of $10^{-8}$), we consistently observe instability of Branch 2 at $\kappa=6.5$; therefore, we have marked a window of instability for Branch 2 for smaller values of $\kappa$ (dashed portion), and marked the branch as stable for larger values of $\kappa$ (solid portion). It is possible, however, that all of Branch 2 is locally stable with a basin of attraction that is shrinking as $\kappa$ approaches the turning point.

A natural and interesting direction is to numerically compute the spectrum of the linearised operator of the PDE about a stationary state, and compare these results with the stability conclusions drawn from the time-dependent solver; note, however, that this requires an appropriate discretisation of the PDE itself, which we avoid by using the map $\mathcal{T}$. We note, however, that this is possible, see for example, \cite{kalise2025linearizationbasedfeedbackstabilizationmckeanvlasov}. Another interesting question is whether there is any relation between the spectrum of the linearised PDE operator and the spectrum of the linearisation of $I - \mathcal{T}$. We leave further explorations of the local and global stability of these branches for future work.

\subsection{The role of the initial iterate}\label{sec:init_iterate}

Newton's method relies heavily on the initial guess provided to the algorithm. To explore the basins of attraction of Algorithm \ref{algorithm:jfnk-fixed-points_intro}, we parametrise initial guesses via two continuous parameters $k,b$:
\begin{align}\label{eq:initial-guesses-suite}
    f_{k,b}(x) = \frac{1}{L} + b\cos\left(\frac{2\pi k x}{L}\right).
\end{align}
This form arises naturally from the expected shape of the locally bifurcating solution given by \eqref{eq:nontrivial-two-mode-solution}, and this is precisely the suite of initial guesses used to obtain all bifurcation diagrams presented thus far.

\begin{figure}[!htbp]
\setlength{\belowcaptionskip}{0pt}
    \centering
    \includegraphics[width=0.75\linewidth]{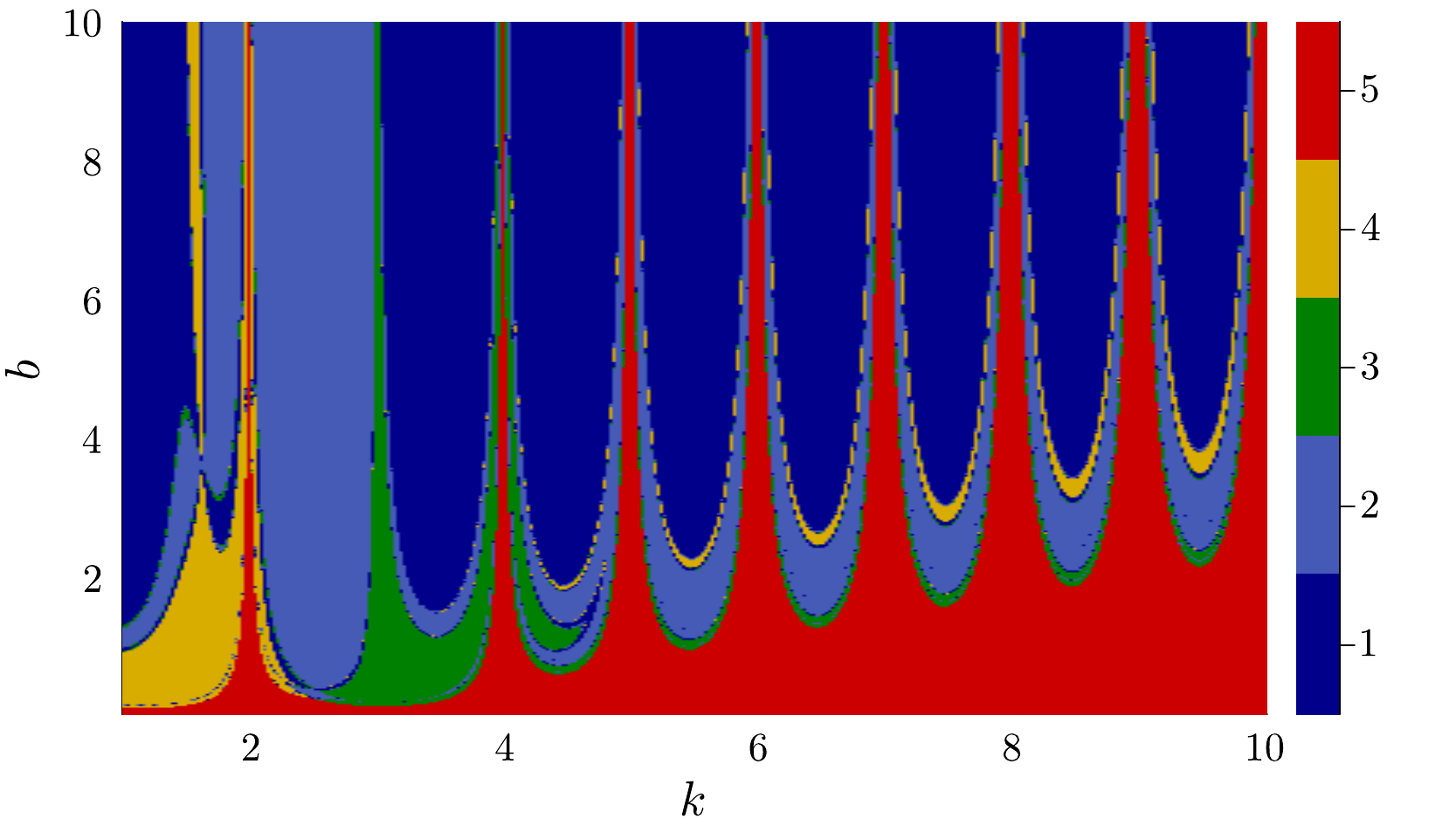}
    \caption{The basins of attraction for initial guesses of the form \eqref{eq:initial-guesses-suite} when fed into Algorithm \ref{algorithm:jfnk-fixed-points_intro}. Note carefully that these are \textbf{not} related to the temporal dynamics of problem \eqref{eq:time_pde_general}. The suite of initial guesses is parametrised by the wavenumber $k$ and the amplitude $b$ as defined in \eqref{eq:initial-guesses-suite}. The heat map indicates which branch the corresponding initial iterate converged to, where the basin colours correspond directly with those branches displayed in Figure \ref{fig:supercritical_bifurcation} (e.g., red regions correspond to those initial guesses which yield Branch 5, the homogeneous solution, in Figure \ref{fig:supercritical_bifurcation}).}
    \label{fig:basins}
\end{figure}

In Figure \ref{fig:basins}, we display the basins of attraction computed for the kernel employed to obtain the supercritical bifurcation diagram from \autoref{fig:supercritical_bifurcation}, at $\kappa = 3$. The colours in Figure \ref{fig:basins} correspond to the solution the algorithm converges to, given that specific pair $(k,b)$. For any integer value $k \not \in \{1,3\}$, we always recover the homogeneous solution. When $k=3$, we recover the three-mode solution predicted by the local bifurcation analysis (provided that $b$ is large enough, otherwise the initial guess is too close to the homogeneous solution to recover an inhomogeneous one). When $k=1$, depending on the size of $b$ we can recover three distinct inhomogeneous solutions: a low-amplitude single-mode solution (yellow) predicted by the local bifurcation analysis, in addition to a high-amplitude single-mode solution (dark blue) and a lower-amplitude multi-modal solution (light blue). For non-integer values of $k$, we observe islands isolated between each integer value of $k$ where we recover all observed solutions. 

This example suggests that the choice of initial iterates is rather delicate, and it is difficult to determine in general whether all possible solutions have been identified. We carefully note that we do not claim that these are all possible solutions, but we can report that for any suite of initial guesses we could produce (e.g., random noise/random perturbations of basis functions, sums of basis functions of differing frequencies, random perturbations of the homogeneous state) are unable to produce additional branches not already obtained via guesses of the form \eqref{eq:initial-guesses-suite}.

\subsection{Two spatial dimensions}\label{sec:twodim_MVE}

\begin{figure}[!htbp]
\setlength{\belowcaptionskip}{0pt}
    \centering
    \includegraphics[width=0.65\textwidth]{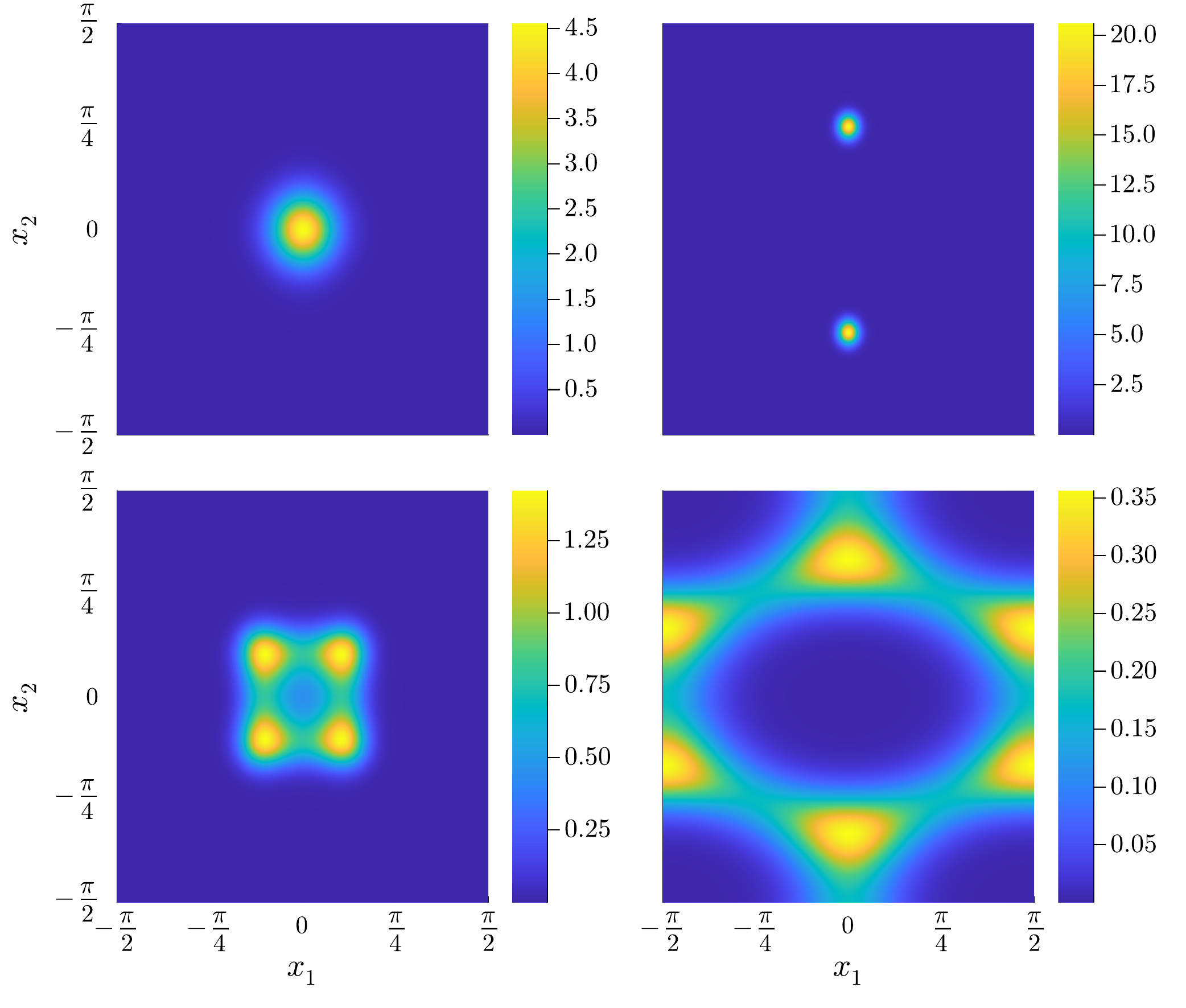}
    \caption{Solutions to the 2-dimensional McKean-Vlasov equation on $\mathbb{T}^2$ of side length $\pi$, each computed in $<6$ seconds. We discretise the domain into $1001 \times 1001$ points. In each case, we chose $\kappa = 10$; the kernel choices were as follows: top left panel: $W(x_1,x_2) = -w_1(x_1) - w_1(x_2)$ (error: $4.09 \times 10^{-14}$); top right panel: $W(x_1,x_2) = -w_1(x_1) - 2w_2(x_1) - w_1(x_2) - 2w_2(x_2)$ (error: $1.16 \times 10^{-10}$); bottom left panel: $W(x_1,x_2) = -w_1(x_1) - 2w_3(x_1) - w_1(x_2) - 2w_3(x_2)$ (error: $5.76 \times 10^{-11}$); and bottom right panel: $W(x_1, x_2) = -\sqrt{2}\, w_1(x_1)w_1(x_2) - \frac{1}{\sqrt{2}}\, w_2(x_1)w_0(x_2)$ (error: $3.74 \times 10^{-11}$).}
    \label{fig:2d_mve}
\end{figure}

Algorithm \ref{algorithm:jfnk-fixed-points_intro} can be naturally adapted to explore solutions of \eqref{eq:McKean_Vlasov_stationary} via \eqref{eq:operator_T_McKean} for dimensions larger than one. In particular, discretisations of domains, functions, and operators can be easily extended across arbitrary dimensions with the caveat that the discretised operator and its Fr{\'e}chet derivative produce a $d$-dimensional tensor, which needs to be flattened into a vector of suitable length before applying the Newton iteration in Algorithm \ref{algorithm:jfnk-fixed-points_intro}. In Figure \ref{fig:2d_mve} we showcase several solution profiles in the two-dimensional case across four distinct kernels. We note that our goal at this stage is to simply demonstrate the applicability in higher dimensions; therefore, we simply report the $\ell^\infty$-error of the discretised $\norm{\mathcal{T}u - u}_\infty$ produced from each of these four stationary profiles. In all cases, the error is below the set convergence tolerance of $10^{-7}$. 

\subsection{The Cucker-Smale and Neural Fokker-Planck Models}\label{sec:CS_NFP_results}

In this subsection, we briefly present our results for the two other problems previously introduced, namely, the Cucker-Smale flocking model \eqref{eq:cucker_smale_stationary} and the neural Fokker-Planck model \eqref{eq:neural_FP_stationary}. For simplicity, we do not provide the same level of technical detail as in the McKean-Vlasov case; the heuristics are identical, using the results of \cite{MR3541988} and \cite{MR4491042, MR4575120}, respectively, to inform our numerical insights. The primary purpose, then, is to emphasise that the approach presented here appears to be very robust to problems defined with different boundary conditions and/or on different domains.

\begin{figure}[!htbp]
\setlength{\belowcaptionskip}{0pt}
  \centering

  \begin{subfigure}[b]{\textwidth}
    \centering
    \includegraphics[width=0.8\linewidth]{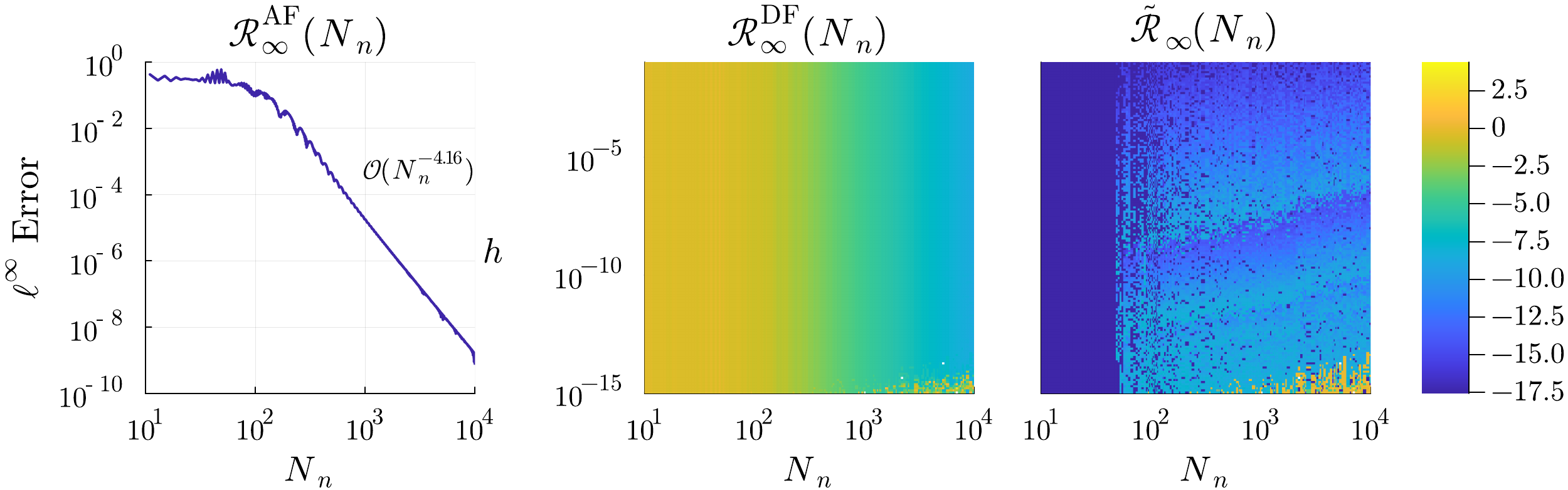}
    \caption{}
    \label{fig:validation_CSM}
  \end{subfigure}

\begin{subfigure}[b]{\textwidth}
    \centering\includegraphics[width=0.8\linewidth]{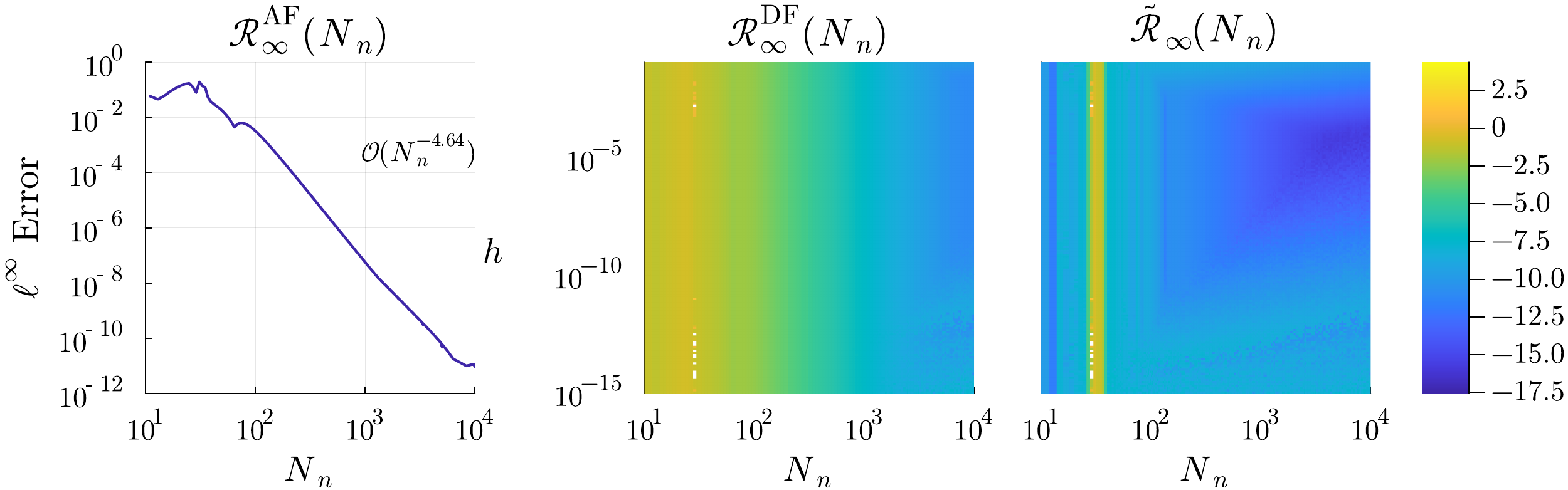}
    \caption{}
    \label{fig:validation_neural_FP}
\end{subfigure}
\caption{Error plots for the numerical solutions obtained to the Cucker-Smale flocking model equation \eqref{eq:cucker_smale_stationary} as described in Section \ref{sec:dev-cucker-smale} (panel (A)) and the homogeneous neural Fokker-Planck equation \eqref{eq:neural_FP_stationary} as described in Section \ref{sec:dev-neural-fp} (panel (B)). For both cases, as an initial guess for Algorithm \ref{algorithm:jfnk-fixed-points_intro}, $f(\cdot) = \exp(-(\cdot)^2/2)$ was used. The values of $N_n, h$ used here are the same as those for Figure  \ref{fig:validation_McKeanVlasov}.}
    \label{fig:validation_CSM_neural_FP}
\end{figure}

First, in Figure \ref{fig:validation_CSM} we emulate Figure \ref{fig:validation_McKeanVlasov} for the Cucker-Smale model \eqref{eq:cucker_smale_stationary}. For this problem, we know there are three possible profiles: a left-skewed solution, a right-skewed solution, and a double-well (symmetric) profile. Therefore, the error computation is always done relative to the corresponding analytical solution (i.e., if we numerically compute the right-skewed solution, error is computed vis-à-vis the analytical right-skewed solution). We have computed errors across initial guesses that favour different solution profiles, and the corresponding error plots are highly similar. We again observe that error decays like $\mathcal{O}(N_n^{-4})$, when using both the analytical and numerical Fr{\'e}chet derivatives.

As further validation, our solver allows us to reproduce \cite[Figure 4]{MR3541988}, which depicts a stability-type region in the $(\kappa, \sigma)$-plane where one expects the first bifurcation to occur (note that in that reference, they use the parameters $(\alpha, D)$). We also applied a crude regression algorithm to predict the power law associated with this curve, and obtained a good estimate of the analytically derived asymptotic limit as the diffusion coefficient vanishes.

Similarly, in Figure \ref{fig:validation_neural_FP} we again emulate Figure \ref{fig:validation_McKeanVlasov} for the neural Fokker-Planck model \eqref{eq:neural_FP_stationary}; the results are essentially identical, with a quartic rate of convergence and good agreement between the use of the analytical versus numerical Fr{\'e}chet derivatives.

\begin{figure}[!htbp]
\setlength{\belowcaptionskip}{0pt}
    \centering
    \includegraphics[width=0.7\linewidth]{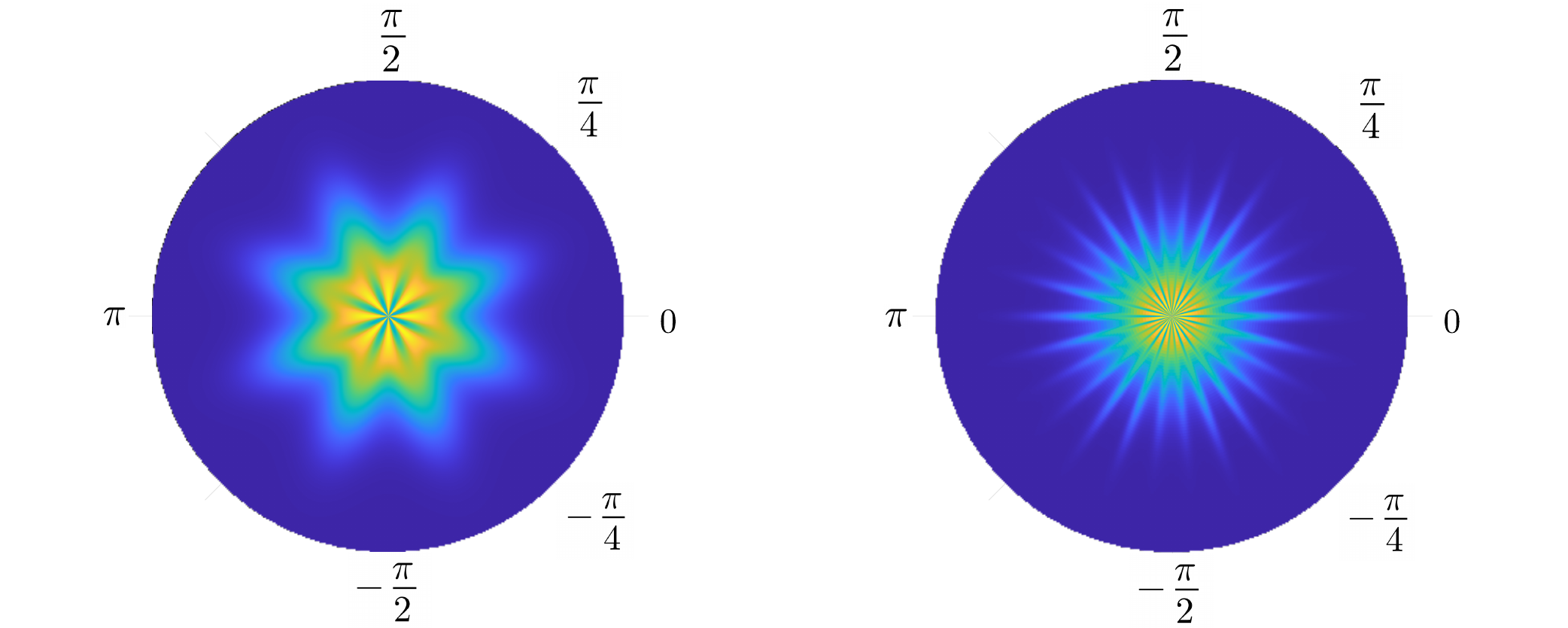}
    \caption{Solution profiles to the neural Fokker-Planck model \eqref{eq:neural_FP_stationary} for two different kernels $W$. We consider $\theta \in (-\pi/2, \pi/2]$ discretised into 1001 samples, and approximate $y \in [0,\infty)$ by considering $y \in [0,35]$, discretised into 3501 samples. We use $B = 3$ and a smoothed ReLU $F(x) = x\max\left(0, \frac{x}{\sqrt{x^2 + 0.1}}\right)$. For the left panel, we use $W(x) = -3 + 3\cos(2x) + 3.3\cos(8x)$ and $\sigma = 0.44$, whereas for the right panel we use a $\tanh$-based kernel $W(x) = -0.005 \cdot 2^{14} \cdot (1 + \tanh(10 - 50|x|))$ defined in \cite{MR4575120}, and $\sigma \approx 0.018805$. These specific values of $\sigma$ have been derived based on the Fourier modes of $W$, as per \cite[Theorem 2.2]{MR4575120}. The initial guesses used to obtain these plots are also based on \cite[Theorem 2.2]{MR4575120}: we define $f_k(x,y) = u_\infty(y) + 0.001\cos\left(\frac{2\pi k x}{L}\right)$. We use $k = 1$ for the left panel and $k = 9$ for the right panel.}
    \label{fig:nfpe_sols}
\end{figure}

Finally, based on the results from \cite{MR4575120}, we display in Figure \ref{fig:nfpe_sols} two inhomogeneous solutions to the neural Fokker-Planck model \eqref{eq:neural_FP_stationary}. Based on the theory of \cite{MR4575120}, the left panel of Figure \ref{fig:nfpe_sols} is consistent with their bifurcation results. In contrast, the right panel of Figure \ref{fig:nfpe_sols} appears to fall outside of the local bifurcation theory and may indicate that a secondary bifurcation has occurred.

\appendix

\subsection*{Acknowledgements} { \footnotesize JAC was partially supported by the Advanced Grant Nonlocal-CPD
(Nonlocal PDEs for Complex Particle Dynamics: Phase Transitions, Patterns and Synchronization) of the European Research Council Executive Agency (ERC) under the European Union Horizon 2020 research and innovation programme (grant agreement No. 883363). JAC acknowledges support by the EPSRC grant EP/V051121/1. YS was partially supported by the Natural Sciences and Engineering Research Council of Canada (NSERC Grant PDF-578181-2023) and graciously acknowledges support from an AARMS Postdoctoral Fellowship. JAC was also partially supported by the “Maria de Maeztu” Excellence Unit IMAG, reference CEX2020-001105-M, funded by MCIN/AEI/10.13039/501100011033/. The authors would also like to thank Lila Spencer, whose summer research project contributed to the implementation of the time-dependent solver used to produce Figure~\ref{fig:triangle_bifurcations}. 
}
\bibliographystyle{abbrv}
\bibliography{references}

\end{document}